\documentclass[secnum,leqno]{amsart}

\usepackage{amsmath}
\usepackage{longtable}
\makeatletter 
\@addtoreset{equation}{subsection} 
\makeatother

\usepackage{mathrsfs, amssymb, array, wasysym}
\usepackage[all,matrix,arrow,curve]{xy}
\usepackage[all]{xy}
\usepackage{mathcomp}

\newcommand{\stackunder}[2]{\mathrel{\mathop{#2}\limits_{#1}}}%
\newcommand{\stackud}[3]{\mathrel{\mathop{#3}\limits_{#1}^{#2}}}%

\newcommand{\type}[1]{$\mathrm{(#1)}$}
\newcommand{\typ}[1]{$\mathrm{#1}$}

\newcommand{\ov}[1]{\overline{#1}}
\newcommand{\ti}[1]{\widetilde{#1}}
\newcommand{\ha}[1]{\widehat{#1}}

\newcommand{\lin}{\text{---}}

\newcommand{\OOO}{\mathscr{O}}
\newcommand{\HHom}{\operatorname{\mathscr{H}\!\mathit{om}}}

\newcommand{\CC}{\mathbb{C}}
\newcommand{\ZZ}{\mathbb{Z}}
\newcommand{\PP}{\mathbb{P}}
\newcommand{\QQ}{\mathbb{Q}}

\newcommand{\toplus}{\mathbin{\tilde\oplus}}
\newcommand{\totimes}{\mathbin{\tilde\otimes}}
\newcommand{\xref}[1]{{\rm \ref{#1}}}

\newcommand{\down}[1]{\left\lfloor #1\right\rfloor}

\newcommand{\Spec}{\operatorname{Spec}}
\newcommand{\Diff}{\operatorname{Diff}}
\newcommand{\wt}{\operatorname{wt}}

\newcommand{\Sing}{\operatorname{Sing}}
\newcommand{\im}{\operatorname{Im}}

\newcommand{\Coker}{\operatorname{Coker}}
\newcommand{\len}{\operatorname{len}}
\newcommand{\Supp}{\operatorname{Supp}}

\newcommand{\red}{\operatorname{red}}
\newcommand{\gr}{\operatorname{gr}}

\newcommand{\ord}{\operatorname{ord}}
\newcommand{\Def}{\operatorname{Def}}
\newcommand{\Hom}{\operatorname{Hom}}
\newcommand{\embdim}{\operatorname{emb\, dim}}

\newcommand{\mt}[1]{\operatorname{#1}}
\newcommand{\muu}{{\boldsymbol{\mu}}}
\newcommand{\compl}{^{\wedge}}
\newcommand{\qq}{\mathbin{\sim_{\scriptscriptstyle{\QQ}}}}

\newtheorem{thm}[subsection]{}

\newtheorem{sthm}[equation]{}

\newtheorem*{lemma*}{Lemma}

\theoremstyle{definition}
\newtheorem{de}[subsection]{}
\newtheorem{sde}[equation]{}
\newtheorem*{remark*}{Remark}

\newcommand{\numer}{\refstepcounter{equation} {\bf \thesthm.}}

\newcounter{eqnumerc}[equation]
\renewcommand{\theeqnumerc}{\rm (\arabic{section}.\arabic{subsection}.\arabic{equation}.\arabic{eqnumerc})}
\newcommand{\eqnumer}{\refstepcounter{eqnumerc} {\theeqnumerc}}

\newcommand{\newpar}[1]{
\par\smallskip\noindent 
\eqnumer\quad
{#1} 
\par\noindent}

\newenvironment{newequation}
{
\begin{equation*}
\refstepcounter{eqnumerc}}
{\leqno{\theeqnumerc}
\end{equation*}
}

\input border.tex
\newcommand{\whitebox}{
\begin{tabular}{c}
\setborder[1cm,0.2cm,0.2cm]\border{\vspace{0.2cm}}
\end{tabular}}

\title{Threefold extremal contractions of type $\text{\type{IA}}$}

\author{
Shigefumi Mori}
\address{RIMS, 
Kyoto University, Oiwake-cho, Kitashirakawa, Sakyo-ku, Kyoto
606-8502, Japan}
\email{ mori@kurims.kyoto-u.ac.jp}
\author{
Yuri Prokhorov}
\address{Department 
 of Algebra, Faculty of Mathematics, Moscow State
 University, Moscow 117234, Russia}
\email{ prokhoro@gmail.com}

\begin{document}

\maketitle
\newcommand\dedicatory[1]{\begin{center}
\texttt{#1}                        
\end{center}
}
\dedicatory{
Dedicated to the memory of Professor Masaki Maruyama}
\bigskip
\begin{abstract} 
Let $(X,C)$ be a germ of a threefold $X$ with terminal singularities
along an irreducible reduced complete curve $C$
with a contraction $f: (X,C)\to (Z,o)$ 
such that $C=f^{-1}(o)_{\red}$ and $-K_X$
is ample. Assume that a general member $F\in |{-}K_X|$
meets $C$ only at one point $P$ and furthermore $(F,P)$ 
is Du Val of type A if $\mathrm{index}(X,P)=4$. 
We classify all such germs in terms of a general member 
$H\in |\OOO_X|$ containing $C$. 
\end{abstract}

\tableofcontents

\section{Introduction}
\begin{de} {\bf Definition.}
\label{setup}
Let $(X,C)$ be a germ of a threefold with terminal singularities
along an reduced complete curve. We say that $(X,C)$
is an \textit{extremal curve germ} if 
there is a contraction $f: (X,C)\to (Z,o)$ such that
$C=f^{-1}(o)_{\red}$ and $-K_X$
is $f$-ample.

If furthermore $f$ is birational, then $(X,C)$ is said to be 
an \textit{extremal neighborhood} \cite{Mori-1988}.
In this case $f$ is called \textit{flipping} if 
its exceptional locus coincides with $C$ 
(and then $(X,C)$ is called \textit{isolated}).
Otherwise the exceptional locus of $f$ is two-dimensional
and $f$ is called \textit{divisorial}.
If $f$ is not birational, then $\dim Z=2$ and 
$(X,C)$ is said to be 
a \textit{$\QQ$-conic bundle germ} \cite{Mori-Prokhorov-2008}.

In this paper, unless explicitly stated otherwise, we assume that $C$ is
\textit{irreducible}.
\end{de}

\begin{de} 
Let $(X, C) $ be an extremal curve germ as above.
For each singular point $P$ of $X$ with $P \in C$, consider
the germ $(P \in C \subset X)$. All such germs (or all such
singular points, for simplicity) are classified into types: 
\type{IA}, 
\type{IC}, \type{IIA}, \type{IIB}, \type{III},
\type{IA^\vee}, 
\type{II^\vee}, 
\type{ID^\vee}, \type{IE^\vee}, as for whose definitions we 
refer the reader to \cite{Mori-1988} and \cite{Mori-Prokhorov-2008}. 
The possible configurations of such points are also classified 
in \cite{Mori-1988} and \cite{Mori-Prokhorov-2008}.
Moreover, it is known that a general member $F\in |{-}K_X|$ has only 
Du Val singularities and all possibilities for $F$ are 
described \cite[7.3, 9.10]{Mori-1988}, \cite[2.2]{Kollar-Mori-1992},
\cite[1.3.7]{Mori-Prokhorov-2008}, \cite[2.1-2.2]{Mori-Prokhorov-2008III}.
The next step in the classification is to study a general hyperplane section, 
that is, a general divisor $H$ of $|\OOO_X|_C$, 
the linear subsystem of $|\OOO_X|$ consisting of sections 
$\supset C$.
Roughly speaking,
the importance of this divisor can be explained as follows. 
Once we have this $H$, the total threefold can be considered as 
a one-parameter deformation of $H$. Then one can apply the deformation theory 
to construct $X$ starting from a two-dimensional data $H\supset C$.

In this paper we classify extremal curve germs of type 
\type{IA} or \type{IA^\vee} in terms of a general member 
$H\in |\OOO_X|_C$. 
An extremal curve germ $(X,C)$ is said to be \textit{of type }
\type{IA} (resp. \type{IA^\vee}) 
if it contains exactly one 
non-Gorenstein point $P$ and it is of type 
\type{IA} (resp. \type{IA^\vee}).
For readers' convenience, we note the 
following characterization (cf. \cite[Th. 2.2]{Kollar-Mori-1992})
for an extremal curve germ $(X,C)$ with a point $P$ of index $m>1$
to be of type \type{IA} or \type{IA^\vee} in terms of a general member 
$F\in |{-}K_X|$: \ $(X,C)$ is of type \type{IA} or \type{IA^\vee}
if and only if (i) ~ $F\cap C=\{P\}$ as a set, and (ii) ~ $(F,P)$ 
is Du Val of type \typ{A} if $m=4$.
\end{de}

\begin{de}
Throughout this paper if we do not specify otherwise 
we assume that $(X,C)$ of type \type{IA} or \type{IA^\vee}. More precisely,
$X$ contains a unique non-Gorenstein terminal point $P\in X$,
which is of type \type{IA} or \type{IA^\vee}. 

A point $(X\supset C\ni P)$ of index $m>1$ is said to be 
of type \type{IA} if there exists an embedding 
$X\subset \CC^4_{x_1,\dots,x_4}/\muu_m(a_1,a_2,-a_1,0)$
such that
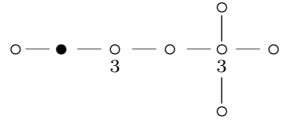
\[
C=\{x_1^{a_2}-x_2^{a_1}=x_3=x_4=0\}/\muu_{m}(a_1,a_2,-a_1,0),
\]
for some positive integers
$a_1$, $a_2$ with $\gcd (a_1a_2, m)=1$ and
$m\in a_1\ZZ_{>0}+a_2\ZZ_{>0}$,
and $X$ is given by an invariant vanishing along $C$ \cite[A.3]{Mori-1988}.
If $f$ is a $\QQ$-conic bundle, then $a_2=1$ by \cite[Proposition 8.5]{Mori-Prokhorov-2008}.
Points of type \type{IA^\vee} are described similarly \cite[A.3]{Mori-1988}.
\end{de}

For a normal surface $S$ and a curve $V\subset S$,
we use the usual notation of graphs $\Delta (S,V)$
of the minimal resolution of $S$ near $V$:
each $\diamond$ corresponds to an irreducible component of $V$ and each
$\circ$ corresponds to an exceptional divisor on the minimal
resolution of $S$, and we may use $\bullet$ instead of $\diamond$ if we want 
to emphasize that it is a complete $(-1)$-curve.
A number attached to a vertex denotes the minus self-intersection number.
For short, we may omit $2$ if the self-intersection is $-2$.

Flipping extremal neighborhoods containing a terminal 
singular point of type \typ{cD/3} (\cite{Mori-1985-cla}, \cite{Reid-YPG1987}) are classified 
in \cite[Theorems 6.2 and 6.3]{Kollar-Mori-1992}.
Thus the following theorem covers the rest of extremal curve germs which 
contain \typ{cD/3} points.

\begin{thm} {\bf Theorem.}
\label{theorem-main-cD/3}
Let $(X,C)$ be an extremal curve germ. Assume that $(X,C)$ is of type \type{IA}
and let $P\in X$ be the non-Gorenstein point.
Assume furthermore that $(X,P)$ is of type \typ{cD/3}.
Then $f$ is a birational contraction, not a $\QQ$-conic bundle.
The general member $H\in |\OOO_X|_C$ and its image $T=f(H)\in |\OOO_Z|$
are normal and have only rational singularities. 
Moreover, if $f$ is not a flipping contraction, then 
the following are the only possibilities for the dual graphs 
of $(H,C)$ and $T$ :

\par\medskip\noindent
$
\begin{array}{lc@{\,}c@{\,}c}
\numer\label{eq-theorem-main-cD/3-A2}\quad \Delta(H,C):
&\circ\lin\bullet\lin\stackrel{3}{\circ}\lin&\circ&\lin\stackrel{3}{\circ}
\\[-4pt]
&&|&
\\[-4pt]
&&\stackunder{3}{\circ}&
\end{array} 
$
\par\noindent
and $T$ is of type \typ{A_2}; here $(X,P)$ is a simple \typ{cD/3} point 
\textup(see \xref{setup-cD/3}\textup);
\par\medskip\noindent
$
\begin{array}{lc@{\,}c@{\,}c}
\numer\label{eq-theorem-main-cD/3-D4}\quad \Delta(H,C):
&&\stackrel{}{\circ}&
\\[-4pt]
&&|&
\\[-4pt]
&\circ\lin\bullet\lin\stackunder{3}{\circ}
\lin\stackunder{}{\circ}\lin&\stackunder{3}\circ
&\lin\stackunder{}{\circ}
\\
&&|&
\\[-4pt]
&&\stackrel{}{\circ}&
\end{array} 
$ 
\par\noindent
and $T$ is of type \typ{D_4}; here $(X,P)$ is a double \typ{cD/3} point.

\par\medskip\noindent
$
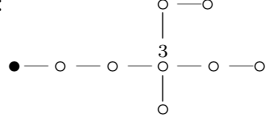
\begin{array}{l@{\,}c@{\,}c@{\,}l}
\numer\label{eq-theorem-main-cD/3-E6}\quad \Delta(H,C):
&&\circ&\lin\circ\\
[-4pt]
&&|\\
&\bullet\lin\circ\lin\circ
\lin&\stackrel{3}{\circ}&\lin\circ\lin\circ\\
[-4pt]
&&|
\\[-4pt]
&&\circ&
\end{array}
$
\par\noindent
and $T$ is of type \typ{E_6}; here $(X,P)$ is a triple \typ{cD/3} point.

\par\medskip
In all the cases above 
the right hand side of the graph for $(H,C)$ corresponds to 
the non-Gorenstein point $P\in H$. The left hand side 
corresponds 
to either a type \type{III} point or a smooth point of $X$. 
\end{thm}
Examples are given in \ref{example-1.1.1} and \ref{example-1.1.2}.


Note that $\QQ$-conic bundles of type \type{IA^\vee} are completely
classified in \cite{Mori-Prokhorov-2008}. The following two theorems 
cover the $\QQ$-conic bundles of type \type{IA}. 

\begin{thm} {\bf Theorem.}
\label{theorem-H-not-normal-conic-bundle}
Let $(X,C\simeq\PP^1)$ be a $\QQ$-conic bundle germ of index $m>2$ and
of type \type{IA}. 
Let $P\in X$ be the non-Gorenstein point.
Then $(X,P)$ is a point of type \typ{cA/m}
and a general member $H\in |\OOO_X|_C$ is not normal.
Furthermore, the dual graph of $(H',C')$, the normalization $H'$
and the inverse image $C'$ of $C$
is of the form
\[
\underbrace{\stackrel{a_r}{\circ}\lin\cdots\lin\stackrel{a_1}{\circ}}_{\varDelta_1}
\lin\bullet\lin 
\underbrace{\stackrel{b_1}{\circ}\lin\cdots\lin\stackrel{b_s}{\circ}}_{\varDelta_2}
\]
\par\noindent
\textup(in particular, $C'$ is irreducible\textup).
Here 
the chain $\varDelta_1$ \textup(resp., $\varDelta_2$\textup) 
corresponds to the singularity of type $\frac1m(1,a)$
\textup(resp., $\frac1m(1,-a)$\textup) 
for some integer $a$ \ \textup($\in [1,m]$\textup) relatively prime to $m$.
The germ $(H,C)$ is analytically isomorphic to the germ 
along the line $y=z=0$
of the hypersurface given by the following weighted polynomial of degree $2m$
in variables $x$, $y$, $z$, $u$:
\[
\phi:= x^{2m-2a}y^2+x^{2a}z^2+yzu.
\]
in 
$\PP(1,a,m-a,m)$. 
Furthermore $(X,C)$ is given as an analytic germ of a subvariety 
of
$\PP(1,a,m-a,m) \times \CC_t$ along $C \times 0$ given by
\[
\phi+ \alpha_1x^{2m-a}y+\alpha_2x^{m-a}uy+\alpha_3x^{2m}+\alpha_4x^{m}u+\alpha_5u^{2}=0
\]
for some $\alpha_1,\ldots \alpha_5 \in t\OOO_{0,{\CC_t}}$ and there is a $\QQ$-conic
bundle structure $X \to \CC^2$ through which the second projection $X \to \CC_t$ factors.
\end{thm}
An explicit example is given in \ref{example-conic-bundle}.

\begin{thm} {\bf Theorem (\cite[\S 3]{Prokhorov-1997_e}, \cite[Th. 12.1]{Mori-Prokhorov-2008}).}
\label{th-index=2}
Let $(X,C\simeq \PP^1)$ be a $\QQ$-conic bundle germ
of index
$2$ and of type \type{IA}. 
Let $f\colon (X,C)\to (Z,o)$ be 
the corresponding contraction. Then $(Z,o)$ is smooth. 
Let $u$, $v$ be a local coordinates on $(Z,o)$. Then there is an embedding
\begin{equation*}
\label{eq-diag-last-2}
f: \xymatrix{X\hspace{4pt} \ar@{^{(}->}[r] 
& {\PP(1,1,1,2)\times Z}
\ar[r]^(0.7){p}
&Z
}
\end{equation*}
such that $X$ is given by two equations
\begin{equation*}
\label{eq-eq-index2}
\begin{array}{l}
q_1(y_1,y_2,y_3)=\psi_1(y_1,\dots,y_4;u,v),
\\[2pt]
q_2(y_1,y_2,y_3)=\psi_2(y_1,\dots,y_4;u,v),
\end{array}
\end{equation*}
where $\psi_i$ and $q_i$ are weighted quadratic in $y_1,\dots,y_4$
with respect to $\wt (y_1,\dots,y_4)=(1,1,1,2)$ and
$\psi_i(y_1,\dots,y_4;0,0)=0$. The only non-Gorenstein point of $X$
is $(0,0,0,1; 0,0)$. Up to projective transformations, the
following are the only possibilities for $q_1$ and $q_2$:

\begin{enumerate}
\item 
$q_1 = y_1^2$, $q_2 = y_2^2-y_1y_3$;\quad
here a general member $H\in |\OOO_X|_C$ is normal.
\item
$q_1 = y_1^2$, $q_2 = y_2^2$; \quad
here every member $H\in |\OOO_X|_C$ is non-normal.
\end{enumerate}
In both cases $C$ is given by $u=v=y_1=y_2=0$.
\end{thm}
Explicit examples are given in Sect. \ref{sect-index-2}
(see also Remark \ref{remark-index-2}).

\begin{de}
\label{explanations-cA/m-cD/3}
For a triple $(X,C,P)$ of type \type{IA} or \type{IA^\vee}, the singularity
$(X,P)$ is either \typ{cA/m}, \typ{cD/3} or of index $2$.
Extremal neighborhoods of index $2$ are classified in 
\cite[\S 4]{Kollar-Mori-1992}. The \typ{cD/3} case is covered 
by \cite[\S 6]{Kollar-Mori-1992} and Theorem \ref{theorem-main-cD/3}.
The next theorem covers the remaining case.
\end{de}

\begin{thm} {\bf Theorem (cf. \cite{Tziolas2005a}).}
\label{theorem-main-birational}
Let $(X,C)$ be an extremal neighborhood of type \type{IA} or \type{IA^\vee}.
Let $P\in X$ be the non-Gorenstein point.
Assume furthermore that $(X,P)$ is of type \typ{cA/m}.
Let $F\in |{-}K_X|$ be a general member. 
Then there exists a member $H\in |\OOO_X|_C$ such that the pair $(X,H+F)$
is LC. 
\begin{sthm}
\label{stheorem-H-normal}
If $H$ is normal, then
$H$ has only log terminal singularities of type \typ{T}.
The graph $\Delta(H,C)$ is of the form 
\begin{newequation}
\label{Esq-theorem-H-normal-birational}
\begin{array}{c@{\,}c@{\,}c@{\,}c@{\,}c@{\,}c@{\,}c@{\,}c@{\,}c@{\,}c@{\,}c@{\,}}
\stackrel{c_1}{\circ}&\lin&\stackrel{c_2}{\circ}&
\lin&\cdots&\lin&\stackrel{c_r}{\circ}&
\lin&\cdots&\lin&\stackrel{c_n}{\circ}
\\[-4pt]
&&&&&&|&&&&
\\[-4pt]
&&&&&&\bullet&&&&
\\[-4pt]
&&&&&&|&&&&
\\[-4pt]
{\circ}&\lin& {\circ}&\lin&\cdots&\lin& {\circ}&&&&
\end{array}
\end{newequation}
Here the chain $[c_1,\dots,c_n]$ corresponds to the non-Du Val singularity 
$(H,P)$ of type \typ{T}. The chain of $(-2)$-vertices in the last line corresponds to a Du Val
point $(H,Q)$. It is possible that this chain is empty \textup(i.e., $(H,Q)$ is smooth\textup).
Cases $r=1$ and $r=n$ are also not excluded.
\end{sthm}

\begin{sthm}
\label{stheorem-H-not-normal}
If every member of $|\OOO_X|_C$ is non-normal, then the dual graph of the normalization $(H',C')$ 
is of the form
\begin{newequation}
\label{Esq-theorem-H-non-normal-birational}
\underbrace{\stackrel{a_r}{\circ}\lin\cdots\lin\stackrel{a_1}{\circ}}_{\varDelta_1}
\lin\bullet\lin
\underbrace{\stackrel{c_1}{\circ}\lin\cdots\lin\stackrel{c_l}{\circ}}_{\varDelta_3}
\lin\diamond\lin
\underbrace{\stackrel{b_1}{\circ}\lin\cdots\lin\stackrel{b_s}{\circ}}_{\varDelta_2}
\end{newequation}
\textup(in particular, $C'$ is reducible\textup).
The chain $\varDelta_1$ \textup(resp., $\varDelta_2$\textup) 
corresponds to the singularity of type $\frac1m(1,a)$
\textup(resp., $\frac1m(1,-a)$\textup) for some $a$ with $\gcd(m,a)=1$ 
and the chain $\varDelta_3$ 
corresponds to the point $(H',Q')$, where $Q'=C_1'\cap C_2'$.
The strings $[a_1,\dots,a_r]$ and $[b_1,\dots,b_s]$ are conjugate \textup(cf. \xref{definition-conjugate}\textup).
Moreover, 
\[
\sum (c_i-2)\le 2 \quad\text{and}\quad \ti C_1^2 +\ti C_2^2 +5-\sum(c_i-2)\ge 0,
\]
where $\ti C=\ti C_1 +\ti C_2$ is the proper transform of $C$ on the minimal resolution
$\ti H$.
Both components of $\ti C$ are contracted on the minimal model 
of $\ti H$.
In this case, the triple $(X,C,P)$ is analytically
isomorphic to $(\{\alpha=0\}, \text{$x_1$-axis}, 0)/\muu_m(1,a,-a,0)$,
where $\gcd(m,a)=1$ and $\alpha(x_1,\dots,x_4)=0$ is the equation of 
a terminal \typ{cA/m}-point in $\CC^4/\muu_m(1,a,-a,0)$. 
\textup(In particular, $(X,C)$ is of type \type{IA}\textup).
\end{sthm}

Conversely, for any germ $(H,C\simeq \PP^1)$ of the form 
\xref{stheorem-H-normal} or \xref{stheorem-H-not-normal}
admitting 
a birational contraction $(H,C)\to (T,o)$ there exists 
a threefold birational contraction 
$f: (X,C)\to (Z,o)$ as in \xref{setup} of type \type{IA}
such that $H\in |\OOO_X|_C$.
\end{thm}

\begin{sde} \textbf{Remark.}
Basically this result is proved in \cite{Tziolas2005a}. 
However \cite{Tziolas2005a}
treated only \emph{divisorial} 
contractions that contracts a divisor to a \emph{smooth} curve. 
Under these assumptions the result of \cite{Tziolas2005a}
is much stronger.
\end{sde}

\begin{sde} \textbf{Remark.}
Note that in \xref{theorem-main-birational} $H$ is not 
assumed to be a general element of $|\OOO_X|_C$.
If $H$ is chosen general, then the cases \ref{stheorem-H-normal}
and \ref{stheorem-H-not-normal} cover all the cases under \ref{theorem-main-birational}.
Proposition \xref{lemma-birational-components-C1} gives a criterion for 
a general member of $|\OOO_X|_C$ to 
be non-normal, and
Proposition \xref{lemma-birational-components-C2} gives, under some additional
assumptions, 
a criterion for a given $H$ to be general.
\end{sde}

To check divisoriality one can use the following criterion which is an
immediate consequence of Theorem \ref{theorem-main-Q-Cartier}.

\begin{thm}{\bf Theorem.}
\label{theorem-main-Q-Cartier-i}
Let $f:(X,C\simeq \PP^1)\to (Z,o)$ be a $3$-dimensional birational
extremal curve germ. 
Then $f$ is divisorial if and only if 
$(Z,o)$ is a terminal singularity. 
\end{thm}

One of our technical tools is the  deformation of extremal curve germs.
In particular, the following result shows that for every
extremal curve germ $f : (X,C) \to (Z,o)$ the contraction $f$
deforms with $X$. Combined with Theorem \ref{theorem-main-Q-Cartier-i}, it allows us to
run MMP for every deformation of an extremal curve germ which may
not be $\QQ$-factorial.

\begin{thm}{\bf Theorem (cf. \cite[(11.4)]{Kollar-Mori-1992}, \cite[(6.2)]{Mori-Prokhorov-2008}).}
\label{theorem-main-def-i}
Let $f : (X,C) \to (Z,o)$ be an extremal divisorial 
\textup{(}resp. flipping, $\QQ$-conic bundle\textup{)} curve germ,
where $C$ is not necessarily irreducible.
Let $ {\pi} : \mathcal{X} \to (\CC_\lambda^1,0)$ be a flat deformation of 
$X= \mathcal{X}_0:= {\pi}^{-1}(0)$
over a germ $(\CC_\lambda^1,0)$ with a flat closed subspace $\mathcal{C} \subset
 \mathcal{X}$ such that $C= \mathcal{C}_0$. Then there exist a
flat deformation $\mathcal{Z} \to (\CC_\lambda^1,0)$ and
a proper $\CC_\lambda^1$-morphism $\mathfrak{f} : \mathcal{X} \to \mathcal{Z}$
such that $f= \mathfrak{f}_0$ and 
$\mathfrak{f}_\lambda : ( \mathcal{X}_\lambda, \mathfrak{f}_\lambda^{-1}(o_\lambda)_{\red})
\to (\mathcal{Z}_\lambda,o_\lambda)$ is a divisorial
\textup{(}resp. flipping, $\QQ$-conic bundle\textup{)} extremal curve germ
for every small $\lambda$, where $o_\lambda:=\mathfrak{f}_\lambda(\mathcal{C}_\lambda)$.
\end{thm}

\begin{de} \textbf{Conventions.}
We work over the complex number field $\CC$.
Notations and techniques of \cite{Mori-1988} will be used freely.
In particular, for a terminal singularity
$(X,P)$ the index-one cover is denoted by 
$(X^\sharp,P^\sharp)\to (X,P)$ and for a subvariety 
$V\subset X$ its preimage is denoted by $V^\sharp$.
\end{de}

\section{Preliminaries}
\begin{de}{\bf Some facts about two-dimensional toric singularities.}
\label{notation-strings}

\begin{sde} {\bf Notation.}
A continued fraction
\[
a_1-\cfrac{1}{a_2-\cfrac{1}{\ddots -\cfrac{1}{a_{r}}}} \qquad
(a_1,\dots, a_r\ge 2)
\]
is denoted by $[a_1,\dots, a_r]$
and called a \textit{string}. Write $m/q=[a_1,\dots,a_r]$, where $\gcd(m,q)=1$.
Given $m$ and $q$ this expression is unique.
It is well-known that the minimal resolution of 
the cyclic quotient singularity $\frac 1m(1,q)$ is 
a chain of smooth rational curves whose self-intersection numbers are 
$-a_1,\dots, -a_r$.
\end{sde}

\begin{sde} {\bf Definition.}
\label{definition-conjugate}
We say that a string $[b_1,\dots,b_s]$ is \emph{conjugate} to $[a_1,\dots, a_r]$
if $[b_1,\dots,b_s]=m/(m-q)$. 
\end{sde}

\begin{sthm} {\bf Lemma.}
\label{lemma-conjugate-strings}
\begin{enumerate}
\item 
If the strings $[a_1,\dots, a_r]$ and $[b_1,\dots,b_s]$ 
are conjugate, then ether $a_1=2$ or $b_1=2$.
\item 
The strings $[a_1,\dots, a_r]$ and $[b_1,\dots,b_s]$ with $a_1=2$ and $r>1$
are conjugate if and only if so are 
$[a_2,\dots, a_r]$ and $[b_1-1,\dots,b_s]$.
\item 
The strings $[a_1,\dots, a_r]$ and $[b_1,\dots,b_s]$
are conjugate if and only if so are 
$[a_r,\dots, a_1]$ and $[b_s,\dots,b_1]$.
\end{enumerate}
\end{sthm}
\end{de}

\begin{de}{\bf \typ{T}-singularities.}
\begin{sde} {\bf Definition (\cite{Kollar-ShB-1988}).}
\label{T-singularities-definition}
A normal surface singularity is said to be of \textit{type 
\typ{T}} if it is log terminal and admits a $\QQ$-Gorenstein one-parameter
smoothing.
\end{sde}

\begin{sthm} {\bf Proposition (\cite[Prop. 5.9]{Looijenga-Wahl-1986}, \cite[Prop. 3.10]{Kollar-ShB-1988}).}
\label{T-singularities-characterization}
A surface singularity is of type \typ{T} if and only if 
it is either Du Val or a cyclic quotient of type 
$\frac 1n (a,b)$, where $\gcd(n,a)=\gcd(n,b)=1$
and $(a+b)^2\equiv 0\mod n$. 
\end{sthm}
\noindent
By \ref{notation-strings}
any non-Du Val \typ{T}-singularity 
is represented by some string $[a_1,\dots,a_r]$.
Then we say that $[a_1,\dots,a_r]$ is a \textit{\typ{T}-string},
or a \textit{string of type \typ{T}}.

\begin{sthm} {\bf Proposition (\cite[Prop. 3.11]{Kollar-ShB-1988}).}
\label{T-singularities-algorithm}
\begin{enumerate}
\item 
The strings $[4]$, $[3,3]$, and $[3,2,\cdots,2,3]$
are of type \typ{T}.
\item
If the string $[a_1,\dots,a_r]$ is of type \typ{T}, then so are
$[2,a_1,\dots,a_{r-1},a_r+1]$ and $[a_1+1,a_2,\dots,a_r,2]$.
\item
Every non-Du Val string of type \typ{T} 
can be obtained by starting with one 
described in \textup{(i)} and iterating the steps described in 
\textup{(ii)}.
\end{enumerate}
\end{sthm}

\begin{sthm} {\bf Corollary.}
\label{Corollary-T-singularities}
Let $(X,P)$ be a $\QQ$-Gorenstein isolated threefold singularity
and let $H\subset X$ be a surface such that $H$ is a Cartier divisor.
If the singularity $(H,P)$ is log terminal, then 
$(H,P)$ is a \typ{T}-singularity and the point $(X,P)$ 
is terminal of type \typ{cA/n} or isolated \typ{cDV}. 
\end{sthm}
\begin{proof}
The only thing we have to prove is the last statement.
By the Inversion of Adjunction
\cite[\S 3]{Shokurov-1992-e-ba}, 
\cite[Ch. 16]{Utah}
the pair $(X,H)$ is PLT.
Since $H$ is Cartier and $(X,P)$ is isolated, it is terminal.
Clearly, we may assume that $(H,P)$ is not Du Val.
Let $F\in |{-}K_X|$ be a general member.
Then $F|_H$ is a general member of $|{-}K_H|$.
Since $(H,P)$ is cyclic quotient (by Proposition 
\ref{T-singularities-characterization}), 
$(H,F|_H)$ is LC.
Again by the Inversion of Adjunction the pair 
$(X,H+F)$ is also LC.
But this means that $(F,P)$ is of type \typ{A}
and so $(X,P)$ 
is of type \typ{cA/n}. 
\end{proof}
\end{de}

\begin{de} {\textbf{Two-dimensional contractions.}}
The following fact is easy and well-known (see, e.g., \cite[Lemma 7.1.11]{Prokhorov-2001}).
\begin{sthm} {\bf Lemma.}
\label{lemma-plt-log-conic-bundle}
Let $\upsilon : S\to R\ni o$ 
be a rational curve fibration germ over a smooth curve and let
$C:=\upsilon^{-1}(o)_{\red}$. 
If the pair $(S,C)$ is PLT, then there is an analytic isomorphism
\[
S \simeq (\PP^1\times \CC)/\muu_m(1,a),
\]
where $\gcd(a,m)=1$.
The graph $\Delta(S,C)$ is of the form
\[
\stackrel{a_r}{\circ}\lin\cdots\lin\stackrel{a_1}{\circ}
\lin\bullet\lin
\stackrel{b_1}{\circ}\lin\cdots\lin\stackrel{b_s}{\circ}
\]
where $[a_1,\dots,a_r]$ and $[b_1,\dots,b_r]$ are conjugate strings.
\end{sthm}

\begin{sthm} {\bf Lemma ({\cite[Th. 6.9]{Shokurov-1992-e-ba}},
{\cite[Prop. 12.3.1-2]{Utah}}).}
\label{lemma-connectedness} 
Let $\upsilon : S\to R$ 
be a rational curve fibration germ over a smooth curve
and let $\Delta$ be an effective $\QQ$-divisor on $S$ such that
$K_S+\Delta\equiv 0$ over $R$. 
Assume 
the locus of log canonical singularities $LCS(S,\Delta)$ of 
$(S,\Delta)$ is not connected near a fiber $\upsilon^{-1}(o)$, $o\in R$. 
Then near $\upsilon^{-1}(o)$
the pair $(S,\Delta)$ is PLT and $\down \Delta$ is a disjoint
union of two sections. 
\end{sthm}
\end{de}

\begin{sthm} {\bf Lemma.}
\label{lemma-surfaces-connect}
Let $C$ be a smooth complete curve contained in a normal surface $H$.
Assume that the pair $(H,C)$ is not PLT at some point, say $P\in C$, 
and $(K_H+C)\cdot C<0$. Then
\begin{enumerate}
 \item 
$H$ has at most two singular points on $C$.
 \item 
If $H$ is singular at a point $Q\in C$ and $Q\neq P$,
then the pair $(H,C)$ is PLT at $Q$.
The dual graph $\Delta(H,C)$ for the 
minimal resolution of $(H,C)$ at $Q$ is of the following form:
\[
\stackunder{C}{\bullet} \lin \underbrace{\stackrel{b_1}{\circ} 
\lin\cdots \lin\stackrel{b_r}{\circ}}_{(H,Q)}
\]
If moreover $(H,Q)$ is a Gorenstein point, then it is Du Val. 
\end{enumerate}
\end{sthm}

\begin{proof}
By the Inversion of Adjunction 
\cite[\S 3]{Shokurov-1992-e-ba}, 
\cite[Ch. 16]{Utah} one has $(K_H+C)|_C=K_C+\Diff_C(0)$, where 
$\Diff_C(0)$ is a $\QQ$-divisor with support at $C\cap \Sing (H)$.
Moreover, the multiplicity of $\Diff_C(0)$ at every 
point of $C\cap \Sing (H)$ is at least $1/2$ and its multiplicity
at $P$ is at least $1$. 
Since $\deg \Diff_C(0)\le -\deg K_C=2$ the assertion of (i) follows.
As for (ii), we see that the multiplicity of $\Diff_C(0)$ at $Q$
is less than $1$. 
Again by the Inversion of Adjunction the pair $(H,C)$ is PLT at $Q$.
The rest follows from the classification of surface PLT pairs
(see, e.g., \cite[Ch. 3]{Utah}).
\end{proof}

\begin{thm} {\bf Lemma.}
\label{lemma-H}
Let $(X,C)$ be an extremal curve germ and let
$f: (X,C)\to (Z,o)$ be the corresponding contraction.
Assume that a member $H\in |\OOO_X|_C$ is normal.
If $(X,C)$ is a $\QQ$-conic bundle germ, then $H$ has only rational singularities.
\end{thm}

\begin{proof}
The assertion
follows from the observation that $H\to f(H)$ is a rational curve fibration.
\end{proof}

\begin{thm} {\bf Theorem 
(\cite[7.3]{Mori-1988}, \cite[1.3.7]{Mori-Prokhorov-2008}).}
Let $(X,C)$ be an extremal curve germ of type 
\type{IA} or \type{IA^\vee} and 
let $P\in X$ be the non-Gorenstein point.
Then a general member $F\in |{-}K_X|$ does not contain $C$
and has only Du Val singularity of type \typ{A} at $P$.
\end{thm}

\begin{thm} {\bf Proposition.}
\label{lemma-lc-near-F}
Let $f: (X,C)\to (Z,o)$ be a contraction
from a threefold with only terminal singularities such that 
$C$ is a \textup(not necessarily irreducible\textup) curve and 
$-K_X$ is ample. Let $F\in |{-}K_X|$ be a general member.
Assume that $F\cap C$ is a point $P$ such that
$(F,P)$ is a Du Val singularity of type~ \typ{A}.
Then, for a general member $H\in |\OOO_X|_C$,
the pair $(X, \, F+H)$ is LC. 

If $f$ is birational, then 
so is the pair $(Z,F_Z+T)$, where $F_Z=f(F)\in |{-}K_Z|$ and $T:=f(H)\in |\OOO_Z|$.
In this case, $(T,o)$ is a cyclic quotient singularity. 
\end{thm}

\begin{proof}
First we consider the case where $f$ is birational
(this case was considered in \cite{Tziolas2005a}).
Then $(F_Z,o)\simeq (F,P)$ is a Du Val singularity of type~\typ{A}.
Let $T$ be a general hyperplane section of $(Z,o)$. Then
$T\cap F_Z$ is general hyperplane section of $(F_Z,o)$.
Clearly, $T\cap F_Z=\Gamma_1+\Gamma_2$ for some irreducible curves $\Gamma_i$
and the pair $(F_Z, \Gamma_1+\Gamma_2)$ is LC.
By the Inversion of Adjunction so is the pair $(Z,F_Z+T)$.
Hence $(T,\Gamma_1+\Gamma_2)$ is LC and 
$(T,o)$ is a cyclic quotient singularity 
(see, e.g., \cite[Ch. 3]{Utah}).
Take $H:=f^*T$. Then $K_X+F+H=f^*(K_Z+F_Z+T)$, i.e., the contraction 
$f$ is $K_X+F+H$-crepant. Hence the pair $(X,F+H)$ is LC.

Now consider the case where $Z$ is a surface.
First we claim that $(X, \, F+H)$ is LC near $F$.
Consider the restriction 
$\varphi=f_F: (F,P) \to (Z,o)$. 
Let $\Xi\subset Z\simeq \CC^2$ be the branch divisor of $\varphi$.
By the Hurwitz formula we can write 
$K_{F}=\varphi^*\bigl(K_Z+\frac12 \Xi\bigr)$.
Hence,
\[
K_{F}+H|_F=\varphi^*\Bigl(K_Z+\frac12 \Xi+T\Bigr).
\]
Using this and the Inversion of Adjunction we get the following equivalences:
$(X, \, F+H)$ is LC near $F$ $\Longleftrightarrow$
$(F, H|_F=\varphi^*T)$ is LC $\Longleftrightarrow$
$(Z=\CC^2, \frac 12 \Xi+T)$ is LC.
Thus it is sufficient to show that $(Z, \frac 12 \Xi+T)$ is LC.

Let $\xi(u,v)=0$ be the equation of $\Xi\subset \CC^2$. 
Then $(F,P)$ is given by the equation
$w^2=\xi(u,v)$ in $\CC^3_{u,v,w}$. 
By the classification of Du Val singularities we can choose 
coordinates $u$, $v$ so that 
\[
\xi=u^2+v^{n+1}. 
\]
Take $T:=\{ v-u=0\}$.
Then $\ord_0 \xi(u,v)|_T=2$.
By the Inversion of Adjunction the pair $(Z,T+\frac12 \Xi)$ is LC.

Thus we have shown that $(X,F+H)$ is LC near $F$.
Assume that $(X,F+H)$ is not LC at some point $Q\in C$.
By the above, $Q\notin F$. 
Note that $H$ is smooth outside of $C$ by Bertini's theorem.

Assume that $H$ is not normal. Let $\nu :H'\to H$ be the normalization and
let $C':=\nu^{-1}(C)_{\red}$.
Write 
$$
K_{H'}+\Diff_{H}(F)=\nu^*(K_X+H+F)\sim 0.
$$
Here $\Diff_{H}(F)=C'+\nu^{-1}(F|_{H})$, where $C'=\nu^{-1}(C)$.
By the Inversion of Adjunction $C'$ is reduced
and $(H',C'+\nu^{-1}(F|_{H}))$ is not LC at $\nu^{-1}(Q)$.
Now we can apply Lemma \ref{lemma-connectedness} to 
$(H',C'+\nu^{-1}(F|_{H})-\varepsilon \upsilon^*(o))$.
\end{proof}

\begin{sthm} {\bf Corollary.}
\label{corollary-structure-H}
Under the assumptions of \xref{lemma-lc-near-F},
if $H$ is not normal, then there is an analytic isomorphism
$(H,P)\simeq \{x_1'x_2'=0\}/\muu_m(a,-a,1)$.
\end{sthm}

\begin{proof}
Let $\pi : (X^\sharp,P^\sharp) \to (X,P)$ be the index-one cover
and let $H^\sharp:=\pi^*H$, $F^\sharp:=\pi^*F$.
Then the pair $(X^\sharp,H^\sharp+F^\sharp)$ is LC. 

Assume that $(X,P)$ is not a cyclic quotient singularity.
One can choose a $\muu_m$-equivariant embedding 
$X^\sharp \subset \CC^4_{x_1,\dots,x_4}$ so that
$\wt (x_1,\dots,x_4)\equiv (a,-a,1,0) \mod m$ and 
$X^\sharp$ is given by the equation $x_1x_2=\phi(x_3^m,x_4)$,
where $\ord_0\phi\ge 2$. 
For some hypersurfaces $D=\{\xi=0\}$ and $S=\{\psi=0\}$ 
in $\CC^4_{x_1,\dots,x_4}$
we have $H^\sharp=D\cap X^\sharp$ and $F^\sharp=S\cap X^\sharp$.
By the Inversion of Adjunction the pair $(\CC^4,X^\sharp+D+S)$ 
is LC. On the other hand, by blowing up the origin we get an
exceptional
divisor of discrepancy 
$$
a(E,X^\sharp+D+S)=3-2-\ord_{0}\xi-\ord_{0}\psi\ge -1.
$$
Hence, $\ord_{0}\xi=1$. Since $\xi$ is an $\muu_m$-invariant,
it contains the term $x_4$. 
Thus $\xi=x_4-\xi'$, where $\ord_0\xi'\ge 2$.
Then $H^\sharp$ is given by two equations
$x_1x_2=\phi(x_3^m,\xi')$ and $x_4=\xi'$.
By changing coordinates we get what we need.

Now assume that $(X,P)$ is a cyclic quotient singularity.
Then $X^\sharp\simeq \CC^3$.
Again one can choose a coordinate system $x_1,x_2,x_3$
in $\CC^3$ so that
$\wt (x_1,x_2,x_3)\equiv (a,-a,1) \mod m$.
Let $\xi$ be the equation of $H^\sharp$.
By blowing up the origin we get $\ord_0\xi\le 2$.
On the other hand, $\xi$ is an invariant.
Hence, $\xi$ contains the term $x_1x_2$
(possibly up to permutations of coordinates if $a\equiv \pm 1$).
\end{proof}

\section{Deformations of $3$-dimensional divisorial contractions}
In this section we recall and set up deformation tools to study extremal curve germs.

\begin{thm}{\bf Theorem.}
\label{theorem-main-Q-Cartier}
Let $f:(X,C)\to (Z,o)$ be a $3$-dimensional divisorial
 extremal curve germ, 
where $C$ is not necessarily irreducible, 
and let $E$ be its exceptional locus. 
Then the divisorial part of $E$ is a $\QQ$-Cartier divisor.
If furthermore $C$ is irreducible, then
$E$ is $\QQ$-Cartier and
 $(Z,o)$ is a terminal singularity. 
\end{thm}

\begin{thm}{\bf Theorem (cf. \cite[(11.4)]{Kollar-Mori-1992}, \cite[(6.2)]{Mori-Prokhorov-2008}).}
\label{theorem-main-def}
Let $f : (X,C) \to (Z,o)$ be an extremal divisorial 
\textup{(}resp. flipping, $\QQ$-conic bundle\textup{)} curve germ,
where $C$ is not necessarily irreducible.
Let $ {\pi} : \mathcal{X} \to (\CC_\lambda^1,0)$ be a flat deformation of 
$X= \mathcal{X}_0:= {\pi}^{-1}(0)$
over a germ $(\CC_\lambda^1,0)$ with a flat closed subspace $\mathcal{C} \subset
 \mathcal{X}$ such that $C= \mathcal{C}_0$. Then there exist a
flat deformation $\mathcal{Z} \to (\CC_\lambda^1,0)$ and
a proper $\CC_\lambda^1$-morphism $\mathfrak{f} : \mathcal{X} \to \mathcal{Z}$
such that $f= \mathfrak{f}_0$ and 
$\mathfrak{f}_\lambda : ( \mathcal{X}_\lambda, \mathfrak{f}_\lambda^{-1}(o_\lambda)_{\red})
\to (\mathcal{Z}_\lambda,o_\lambda)$ is a divisorial
\textup{(}resp. flipping, $\QQ$-conic bundle\textup{)} extremal curve germ
for every small $\lambda$, where $o_\lambda:=\mathfrak{f}_\lambda(\mathcal{C}_\lambda)$.
\end{thm}

\begin{sthm}{\bf Corollary.}
Let $f : (X,C) \to (Z,o)$ be an extremal divisorial curve germ,
where $C$ is not necessarily irreducible.
Let $P^{(1)},\dots, P^{(r)}\in X$ 
be singular points. Let
$(X_\lambda,P_\lambda^{(i)})\supset (C_\lambda,P_\lambda^{(i)})$ 
be a set of local $1$-parameter analytic deformations of $(X,P^{(i)})\supset (C,P^{(i)})$.
Then it extends to a $1$-parameter analytic deformation
$X_\lambda\supset C_\lambda\supset \{P_\lambda^{(1)},\ldots, P_\lambda^{(r)}\}$ of global 
$X\supset C\supset \{P^{(1)},\ldots, P^{(r)}\}$ in the sense that there exist a
flat deformation $\mathcal{Z} \to (\CC_\lambda^1,0)$ and
a proper $\CC_\lambda^1$-morphism $ \mathfrak{f} : \mathcal{X} \to \mathcal{Z}$
such that $f= \mathfrak{f}_0$ and 
$\mathfrak{f}_\lambda : ( \mathcal{X}_\lambda, \mathfrak{f}_\lambda^{-1}(o_\lambda)_{\red})
\to (\mathcal{Z}_\lambda,o_\lambda)$ is a divisorial extremal curve germ
for every small $\lambda$, where $o_\lambda:=\mathfrak{f}_\lambda(\mathcal{C}_\lambda)$.
\end{sthm}

We need the following easy lemma which can be found in
\cite[(9.3)]{Bingener1981} (without proof).

\begin{thm}{\bf Lemma.} 
\label{lemma-def-analytic}
Let $p:{\mathcal D} \to \mathcal{X}\supset \ell$ 
be an arbitrary analytic morphism 
and let $\ell \subset X$ be a compact 
subset such that  
$p^{-1}(\ell)$ is compact. Then 
there exist open subsets $W \supset p^{-1}(\ell)$ of ${\mathcal D}$
and $V \supset p(W)$ of $\mathcal{X}$ such that 
$p|_W: W\to V$ is proper and
$p(W) $ is
an analytic subset of $V$.
\end{thm}

\begin{proof}
There is an open subset $U \supset p^{-1}(\ell)$ of ${\mathcal D}$
such that $\bar U$ is compact (and $U$ is open and closed in
${\mathcal D} \setminus \partial U$).
Since $p(\partial U)$ is a closed set disjoint from $\ell$, there is an
open set $V \supset \ell$ such that $\bar V$ is disjoint from $p(\partial U)$.
Then $p^{-1}(\bar V)$ is disjoint from $\partial U$.
Hence $W:=U \cap p^{-1}(V)$ is an open and closed subset of $p^{-1}(V)$ 
and is $\bar W \subset U$ is compact.
Hence $p|_{W} : W \to V$ is proper.
This means that $p(W)$ is an analytic subset of $V$.
\end{proof}

The following is the key step in the proof of \ref{theorem-main-Q-Cartier}
and \ref{theorem-main-def}.

\begin{thm}{\bf Proposition.}
\label{deformation-curve-germ}
Let $f : (X,C) \to (Z,o)$ be a divisorial extremal curve germ,
where $C$ is not necessarily irreducible.
Let $\bar {\pi} : \bar{\mathcal{X}} \to (\CC_\lambda^1,0)$ be a flat deformation of 
$X=\bar{\mathcal{X}}_0:=\bar {\pi}^{-1}(0)$
over a germ $(\CC_\lambda^1,0)$.

\textup{(i)}
Let $\bar{\mathcal{X}}\compl$ be the completion of 
$\bar{\mathcal{X}}$ along $\lambda=0$. Then $f: X\to Z$ extends to a contraction 
$\mathfrak{f}\compl: \bar{\mathcal{X}}\compl  \to {\mathcal{Z}}\compl$.

\textup{(ii)} 
Let $n$ be an arbitrary positive integer.
Then there exist flat deformations ${\pi} : \mathcal{X} \to (\CC_\lambda^1,0)$ 
and $\mathcal{Z} \to (\CC_\lambda^1,0)$ and 
a proper $\CC_\lambda^1$-morphism 
$\mathfrak{f} : \mathcal{X} \to \mathcal{Z}$
such that ${\pi}_{(n)} \simeq \bar {\pi}_{(n)}$,
$f=\mathfrak{f}_0$, and 
$\mathfrak{f}_\lambda : \mathcal{X}_\lambda
\to \mathcal{Z}_\lambda$
is a divisorial contraction
\textup(which contracts a divisor to a curve\textup) 
for every small $\lambda$, 
where 
$\mathcal{A}_{(i)}:=\mathcal{A} \times_{\CC_\lambda^1} 
\operatorname{Spec} \CC [[\lambda]]/(\lambda^{i+1})$ for any
object $\mathcal{A}$ over $\CC_\lambda^1$ and $i \ge 0$.
\end{thm}

\begin{proof}
Let $\phi \in H^0(X,\OOO_X)$ be
a general section vanishing on $C$ and $H$ (resp. $H_Z$) the member of
$|\OOO_X|$ (resp. $|\OOO_Z|$) defined by $\phi$ (resp. $f_*\phi$).
We note that $H$ (resp. $H_Z$) is smooth outside $C$ (resp. $o$)
and $f$ induces an isomorphism 
$H \setminus C \simeq H_Z \setminus \{o\}$.

Then as in \cite[(11.3)--(11.4)]{Kollar-Mori-1992}, the miniversal deformation 
spaces $\Def(H)$ and 
$\Def(H_Z)$ exist as analytic spaces and $f$ induces 
a complex analytic morphism $\Def(f,H) : \Def(H) \to \Def(H_Z)$.
Let $\phi: X\to \CC_s^1$ be the morphism defined by 
$s=\phi$. This morphism is a flat family of $H$ over 
$\CC_s^1$. Thus we have an induced morphism 
$\bar w: \CC_s^1\to \Def (H)$, i.e. an element
$\bar w\in \Hom(\CC_s^1, \Def (H))$.
Furthermore $X$, $Z$, and $f$ can be reconstructed by 
the morphism
$\bar w: (\CC_s^1,0) \to \Def(H)$.
Our goal is to construct the following morphism 
extending $\bar w$.
\begin{equation*}
\label{q-diagram-cD-2}
w: (\CC_{s,\lambda}^2,0) \longrightarrow\Def(H)
\end{equation*}

Since $R^1f_*\OOO_X=0$, the section $\phi$ extends to a
formal  section $\hat{\phi}$ on
the completion $\bar{\mathcal{X}}\compl$ of $\bar{\mathcal{X}}$ along $X$.
This proves (i).
We thus see that 
$\bar w\in \Hom(\CC_s^1,\Def(H))$ 
extends to
$\hat w\in \Hom((\CC_{s,\lambda}^2,0)\compl,\Def(H))$, where
$({\CC_{s,\lambda}^2},0)\compl$ is the completion of $(\CC_{s,\lambda}^2,0)$ along $\{\lambda=0\}$.
Then by \cite[Theorem 1.5, (i)]{Artin1968}, 
$\hat w$ can be approximated by
an analytic extension   
$w\in \Hom((\CC_{s,\lambda}^2,0),\Def(H))$ 
of $\bar w$.
This gives us a flat family 
$\mathcal{X}$ 
over $\CC^1_\lambda$ approximating 
$\bar{\mathcal{X}}$.

It remains to settle divisoriality. Arbitrarily close to $C$ 
there is an $f$-exceptional curve $\ell\simeq\PP^1$ 
such that $N_{\ell/X} \simeq \OOO_{\ell} \oplus \OOO_{\ell}(-1)$ which sweep out
an $f$-exceptional divisor of $X$. Hence 
$N_{\ell/\mathcal{X}} \simeq \OOO_{\ell}^{\oplus 2} \oplus \OOO_{\ell}(-1)$ and
there are no obstructions to deforming these $\ell$ out to
$\mathcal{X}_\lambda$ and hence $\mathfrak{f}_\lambda$ contracts a divisor.
This proves the statement (ii) of our proposition.
\end{proof}

\begin{proof}[Proof of Theorem \xref{theorem-main-Q-Cartier}]
Let $P^{(1)},\dots, P^{(r)}\in X$ 
be singular points.
As in \cite[Appendix 1b]{Mori-1988}
one can see that 
every local deformation of singularities extends 
to a deformation  of global $X$. 
For every terminal singularity $(X,P^{(i)})$ we take a 
$\QQ$-smoothing, a deformation whose general member
has only cyclic quotient singularities \cite[(6.4)]{Reid-YPG1987}.  
By the above, there exists a one-parameter deformation 
$\bar{\mathcal{X}}$ over a disk in $\CC^1_\lambda$ such that  $\bar{\mathcal{X}}_0\simeq X$
and, for small $\lambda\neq 0$, the fiber $\bar{\mathcal{X}}_{\lambda}$ has only
terminal cyclic quotient singularities.
Then we apply Proposition \ref{deformation-curve-germ}, (ii).
In notation of Proposition \ref{deformation-curve-germ},
there exists a divisorial contraction $\mathfrak{f}: \mathcal{X}\to \mathcal{Z}$ 
contracting a divisor $\mathcal{E}$
(the divisorial part of the exceptional locus) 
to a surface on $\mathcal{Z}$ and, 
for small $\lambda\neq 0$, the fiber ${\mathcal{X}}_{\lambda}$ also has only
terminal cyclic quotient singularities 
because at every singular point $P$ of $X$ the local germ of $\bar{\mathcal{X}}$ at $P$
can be approximated by one of $\mathcal{X}$ to an arbitrarily high order of $\lambda$.

Let $P\in X= \mathcal{X}_0$ be a singular point and 
let $(X^\sharp, P^\sharp)$ be the index-one cover.
Then the local deformation  $(\mathcal{X},P)$ 
is induced by a deformation $(\mathcal{X}^\sharp,P^\sharp)$ of
$(X^\sharp, P^\sharp)$ (cf. \cite[The last paragraph of \S 6]{Stevens-1988}). 
Since the germ $(X^\sharp, P^\sharp)$ is a
hypersurface singularity  \cite{Reid1983}, so is 
$(\mathcal{X}^{\sharp},P^\sharp)$. 
Moreover, the singularity $(\mathcal{X}^{\sharp},P^\sharp)$
is isolated. 
Hence by \cite[Exp. XI, Corollary 3.14]{Grothendieck1968} the variety 
$\mathcal{X}^\sharp$ is factorial at $P^\sharp$ and so
$\mathcal{X}$ is $\QQ$-factorial at $P$. In particular, 
$\mathcal{E}$ is a $\QQ$-Cartier divisor.
Thus $\mathcal{E}|_X=E$ on $X\setminus C$.
If moreover $C$ is irreducible, then $\rho(X)=1$ (see \cite[(1.3)]{Mori-1988}) 
and so $K_X\qq\mathcal{E}|_X$.
Hence, $\mathcal{E}|_X$ is negative on $C$ and $\mathcal{E}|_X\supset C$.
This implies that $E=\mathcal{E}|_X$ and it is also $\QQ$-Cartier.
\end{proof}

\begin{proof}[Proof of Theorem \textup{\ref{theorem-main-def}}]
The flipping case follows from \cite[(11.4)]{Kollar-Mori-1992} and
the $\QQ$-conic bundle case from \cite[(6.2)]{Mori-Prokhorov-2008}.
So, we assume that $f$ is divisorial. 
Let $E \subset X$ be the exceptional divisor of $f$ and let 
$E_i$'s be its irreducible components.
Then, for each $i$, $B_i:=f(E_i) \subset Z$ is an 
irreducible curve passing through $o$.

First, we treat th case where $C$ is irreducible.
Then by Theorem \ref{theorem-main-Q-Cartier} \ $E$ is a $\QQ$-Cartier 
divisor and $Z\ni o$ is a terminal singularity.

For each $E_i$, choose a smooth fiber $\ell'_i$ of $E_i \to B_i$
and let $[\ell'_i]$ degenerate to $[\ell_i]$ lying over $o$ 
in the Douady space of $X/Z$. We assume each $[\ell'_i]$
is chosen arbitrarily close to $[\ell_i]$.
Consider the closed subspace $A'$ of the Douady space of $X/Z$
parametrizing all compact subspaces $F \subset X$ with
$\Supp F \subset C$. Then each irreducible component of $A'$ 
is compact, cf. \cite{Fujiki1978/79}, and let $A$ be the smallest open and closed subset
of $A'$ containing all $[\ell_i]$. Thus $A$ is also compact.
Then we work on a sufficiently small neighborhood ${\mathcal D'}$ of 
$A$ in the Douady space of $\mathcal{X}/\CC^1_{\lambda}$
such that ${\mathcal D'} \ni [\ell'_i]$ for each $i$.
 
We note that $\mathcal{X}$ is smooth along each $\ell'_i$ and
$N_{\ell'_i/\mathcal{X}} \simeq \OOO_{\ell'_i}^{\oplus 2} \oplus
\OOO_{\ell'_i}(-1)$.
Hence ${\mathcal D}'$ is smooth of dimension $2$ at each
$[\ell'_i]$. 
Let ${\mathcal D} \subset {\mathcal D}'$ be the smallest one
among the union of the irreducible components of $\mathcal{D}'$ 
such that
${\mathcal D}\ni [\ell'_i]$ for all $i$.
Then ${\mathcal D}$ is a two-dimensional
closed subspace of ${\mathcal D}'$.

Let $\mathcal{T} \subset \mathcal{X} \times_{\CC_{\lambda}^1} {\mathcal D}$ 
be the universal closed subspace parametrized by ${\mathcal D}$ 
with two projections
$\pi : \mathcal{T} \to {\mathcal D}$ and $p : \mathcal{T} \to \mathcal{X}$.

We note that $p^{-1}(C)\subset A$ and it is  compact because
the variety
$\mathcal{X}_0=X$ has a divisorial contraction to $Z$,
 $C$ is the fiber over $o \in Z$,
and $\pi^{-1}(t)$ does not intersect $C$ for $t\notin A$.

Let $\mathcal{E} :=p(\mathcal{T}) \subset \mathcal{X}$ be the image
of the proper morphism $p$ and
it is an analytic subset by Lemma \ref{lemma-def-analytic}.
We also denote by $p : \mathcal{T} \to \mathcal{E}$ the morphism
induced by $p$ and let 
$p : \mathcal{T} \stackrel{p'}{\longrightarrow} 
\bar {\mathcal{T}}\stackrel{\bar p}{\longrightarrow}\mathcal{E}$
be
the Stein factorization of $p$ so that $p'_*\OOO_{\mathcal{T}}=\OOO_{\bar {\mathcal{T}}}$.

\begin{sthm}{\bf Claim.}
$\mathcal{E}$ is a $\QQ$-Cartier divisor.
\end{sthm}

\begin{proof}
Let $\mathcal{X}\compl$ be the completion of $\mathcal{X}$
along $\lambda=0$.
By Proposition \ref{deformation-curve-germ}, (i) the morphism
$f: X\to Z$ extends to a contraction 
$\mathfrak{f}\compl: {\mathcal{X}}\compl  \to {\mathcal{Z}}\compl$,
where $\mathcal{Z}\compl$ is $\QQ$-Gorenstein
\cite[Corollary 10]{Stevens-1988}
because $Z$ is terminal.
Comparing $K_{{\mathcal{X}}\compl}$ and ${\mathfrak{f}\compl}^*K_{{\mathcal{Z}}\compl}$
we see that there is an effective $\QQ$-Cartier divisor 
$\mathcal{F}\compl\qq K_{{\mathcal{X}}\compl}-
{\mathfrak{f}\compl}^*K_{{\mathcal{Z}}\compl}$ on ${\mathcal{X}}\compl$
such that $\mathcal{F}\compl|_{{X}\compl}={E}\compl$ and
$\mathcal{F}\compl=\mathcal{E}\compl$ outside of ${C}\compl$.
Hence $\mathcal{F}\compl=\mathcal{E}\compl$.
\end{proof}

Now we define a morphism $q : \mathcal D \to \mathcal{B}$ 
such that $q(p^{-1}(C))$ is one point as follows.
Take a general point $\zeta$ of $C$ and take a small $3$-dimensional
disk $(\Delta^3,0)$ centered at $\zeta$ and transversal to $C$ at $\zeta$.
Then the Cartier divisor $\Delta^3$ in a neighborhood of $C$ induces
a Cartier divisor of $\mathcal T$ finite and flat over $\mathcal D$.
Let $d$ be the degree of $p^{-1}(\Delta^3)/\mathcal D$.
Then $x \in \mathcal D \mapsto \pi^{-1}(x)\cap p^{-1}(\Delta^3)$
associates to $x$ a $0$-cycle of degree $d$ on $\Delta^3$
and we have thus a required morphism 
$q : \mathcal D \to \mathcal{B} :=S^d(\Delta^3)$
such that $q(p^{-1}(C))$ is the $0$-cycle $d \cdot [0]$.

We claim that we have a proper morphism $r : \bar {\mathcal{T}} \to \mathcal{B}$
making the following diagram commutative
\[
\xymatrix{
\mathcal T \ar[r]^{p'}\ar[d]^{\pi} &{\bar {\mathcal{T}}}\ar[d]^r \ar[r]^{\bar p}&\mathcal{E}
\\
\mathcal D\ar[r]^{q} & \mathcal{B} &
}
\]
Indeed since $q(\pi(p^{-1}(C)))$ is one-point $d \cdot [0]$, 
we can shrink $\mathcal E$
so that $q(\mathcal D)$ is contained in a Stein open neighborhood
of $d \cdot [0]$. Hence the morphism $\mathcal T \to \mathcal{B}$
factors through $p' : \mathcal T \to \bar{\mathcal{T}}$ and the
claim is proved.

We claim that $p$, $p'$, $\bar p$ are isomorphisms over every $\ell'_i$
and in particular $\bar p$ is finite and bimeromorphic.
Indeed by $N_{\ell'_i/\mathcal{X}} \simeq \OOO_{\ell'_i}^{\oplus 2} \oplus
\OOO_{\ell'_i}(-1)$, $p$ is an isomorphism near $\pi^{-1}([\ell'_i])$
and by the divisorial contraction on $X=\{\lambda=0\} \subset \mathcal X$
one has $p^{-1}(\ell'_i)=\pi^{-1}([\ell'_i])$. These settle the claim.

Let $\mathfrak{c}:=\HHom_{\OOO_{\mathcal{T}}}(\OOO_{\bar{\mathcal{T}}},
\OOO_{\mathcal{T}})$
be the conductor of $\bar p$ and let $V(\mathfrak{c}) \subset \bar {\mathcal{T}}$
be the locus defined by $\mathfrak{c}$. Then we claim that 
$r(V(\mathfrak{c}))$ is finite over $\CC^1_{\lambda}$.
Indeed this is obvious since 
$r(V(\mathfrak{c}))\not\ni q([\ell_i'])$ and 
the fiber of $r(V(\mathfrak{c}))$ over $\{\lambda=0\}$
is a finite set.

Let $J \subset \OOO_{\mathcal{B}}$ be an arbitrary sheaf of ideals such that
$J\OOO_{\bar{\mathcal{T}}} \subset \mathfrak{c}$
and  $V(J)$ is finite over $\CC_{\lambda}^1$.
By 
\cite[Theorem (6.1)]{Bingener1981} 
we have the following digram 
\[
\xymatrix{
V(J)\ar[r]\ar@{_{(}->}[d]&\CC^1_{\lambda}\ar[d]
\\
\mathcal{B}\ar[r]^{q'}&\mathcal{E}'
} 
\]
where $\mathcal{E}':=\mathcal{B} \coprod_{V(J)}\CC^1_{\lambda}$
is the amalgamated sum
(coproduct) of $\mathcal{B}$ and $\CC^1_{\lambda}$ over $V(J)$ and 
$q' : \mathcal{B} \to \mathcal{E}'$ is  a bimeromorphic finite morphism.
Since $\mathfrak{c}$ is the conductor of $\bar p$, we have  
\[
\mathcal{E}=\bar{\mathcal{T}} \coprod\nolimits_{V_{\bar{\mathcal{T}}}(\mathfrak{c})}
V_{\mathcal{E}}(\mathfrak{c})
\]
 and the following commutative diagram.
\[
\xymatrix{
V_{\bar{\mathcal{T}}}(\mathfrak{c})\ar[r]\ar@{_{(}->}[d]
& V_{\mathcal{E}}(\mathfrak{c}) \ar@{_{(}->}[d]
\\
\bar{\mathcal{T}}\ar[r]^{\bar p}&\mathcal{E}
} 
\]
These two diagrams fit into a big one which allows us to define 
an induced morphism $\eta:  \mathcal{E}\to \mathcal{E}'$:
\[
\xymatrix{
V_{\bar{\mathcal{T}}}(\mathfrak{c})\ar[r]\ar@{_{(}->}[d]\ar[drr]
& V_{\mathcal{E}}(\mathfrak{c}) \ar@{_{(}->}[d]\ar@/^/[drr]^{r}
\\
\bar{\mathcal{T}}\ar[r]\ar@/_/[drr]
&\mathcal{E}\ar@{-->}@/^4pt/[drr]_<(0.3){\eta}&
V(J)\ar[r]\ar@{_{(}->}[d]
&\CC^1_{\lambda}\ar[d]
\\
&&\mathcal{B}\ar[r]_{q'}
&\mathcal{E}'
} 
\]
Finally we have the following commutative diagram:
\[
\xymatrix{
\mathcal T\ar[r]^{p'}\ar[d]^{\pi} 
&{\bar {\mathcal{T}}}\ar[r]^{\bar p}\ar[d]^{r} 
&\mathcal{E}\ar[d]^{\eta}
\\
\mathcal D\ar[r]^{q} 
&\mathcal{B}\ar[r]^{q'} & \mathcal{E}'
}
\]

For any $i,\, j\ge 0$ 
the sheaf $\OOO_{i\mathcal{E}}(-j\mathcal{E})$ denotes  the quotient
$\OOO_{\mathcal{X}} (-j\mathcal{E}) /\OOO_{\mathcal{X}}(-(i+j)\mathcal{E})$.
\begin{sthm}{\bf Claim.}
\label{claim-new-deformation}
For any $i,\, j\ge 0$ we have $R^1\eta_*\OOO_{i\mathcal{E}}(-j\mathcal{E})=0$.
Therefore, the following sequence 
\[
 0\longrightarrow
\eta_*\OOO_{i\mathcal{E}}(-j\mathcal{E})
\stackrel{}{\longrightarrow}
\eta_*\OOO_{(i+j)\mathcal{E}}
\longrightarrow
\eta_*\OOO_{j{\mathcal{E}}}
\longrightarrow 0.
\]
is exact. 
\end{sthm}
\begin{proof}
By the Kawamata-Viehweg vanishing \cite[Theorem 3.6]{Nakayama1987}
$R^1f_*\OOO_X(-kE)=0$ for $k\ge 0$.
Then from the exact sequence
\[
0\longrightarrow \OOO_{{X}}(-(i+j){E})
 \longrightarrow \OOO_{{X}}(-j{E})
\longrightarrow \OOO_{iE}(-j{E})
\longrightarrow 0
\]
we see that $R^1f_*\OOO_{iE}(-j{E})=0$ for $i,\, j\ge 0$.
Now we assert that the following sequence
\[
0\longrightarrow
\OOO_{i\mathcal{E}}(-j\mathcal{E})
\stackrel{\cdot \lambda}{\longrightarrow}
\OOO_{i\mathcal{E}}(-j\mathcal{E})
\longrightarrow
\OOO_{iE}(-j{E})
\longrightarrow 0.
\]
is exact for $i,\, j\ge 0$.
Recall that the space $\mathcal{X}$ is $\QQ$-Gorenstein
\cite[The last paragraph of \S 6]{Stevens-1988}.
Consider 
the index-one cover 
$\nu: (\mathcal{X}^\sharp, P^\sharp)
\to (\mathcal{X}, P)$ 
with respect to $\mathcal{E}$ at an arbitrary point $P\in X$.  
Since the map $\nu$ is \'etale in codimension two, 
both $\mathcal{X}^\sharp$ and 
$X^\sharp: =\nu^{-1}(X)$ are terminal.
The induced divisors $\mathcal{E}^\sharp$ 
and $E^\sharp$ are Cartier on $\mathcal{X}^\sharp$ 
and $X^\sharp$, respectively, and $E^\sharp= \mathcal{E}^\sharp|_{X^\sharp}$.
Hence the assertion on exactness can be readily checked on $\mathcal{X}^\sharp$. 
Then by Nakayama's lemma we obtain 
$R^1 \eta_*\OOO_{i\mathcal{E}}(-j\mathcal{E})=0$.
\end{proof}

Fix a positive integer $m$ such that 
 both $mE$  and $m\mathcal{E}$ are Cartier and
define a ringed space $\mathcal{E}''$ as a topological space 
$\operatorname{Spec}_{\mathcal{E'}}\eta_*\OOO_{\mathcal{E}}$
with the sheaf of rings $\eta_*\OOO_{m\mathcal{E}}$. 
Then $\mathcal{E}''$ is a complex space by Claim \ref{claim-new-deformation}
and \cite[\S 10]{Bingener1981}.

Now we would like to show that  $\mathcal X$ 
has a modification, and in order to do that we would like to check the 
conditions (1) and (2) of 
Corollary (8.2)
of Bingener \cite{Bingener1981} 
for the morphism 
$\mathcal{X}\supset  m\mathcal{E} \to 
\mathcal{E}''$
induced by $\eta$.
The condition (1) is obvious because $-\mathcal{E}$ is ample,
and the condition  (2) follows from 
the exact sequence in
Claim \ref{claim-new-deformation} with $j=1$.
Thus the desired contraction $\mathfrak{f}: \mathcal{X}\to\mathcal {Z}$
exists by \cite[Corollary (8.2)]{Bingener1981}.
So the proof of the case of irreducible $C$ is completed.

Now we consider the general case, i.e. we assume that 
$C$ is reducible.
Run an analytic MMP on $\mathcal{X}$ in the following way.
Every irreducible $K$-negative curve on the central fiber 
of $X/Z$ generates an extremal ray on $X$.
By  \cite[(11.7)]{Kollar-Mori-1992} flips on $X$ extend to ones on $\mathcal{X}$.
So do divisorial contractions by our previous arguments.
By Theorem \ref{theorem-main-Q-Cartier} we stay in the terminal category.
At the end we get $X'\subset \mathcal{X}'/\CC_\lambda^1$
such that $X'$ is a minimal model over $Z$.
Moreover, all fibers of $f': X'\to Z$ are of dimension $\le 1$ and
$-K_{X'}$ is ample over $Z$ outside of the central fiber. Hence,  
$f': X'\to Z$ is a small contraction.
Note that $R^1f'_*\OOO_{X'}=R^1f_*\OOO_{X}=0$.
By \cite[(11.4)]{Kollar-Mori-1992}  the contraction $f': X'\to Z$
extends to $\mathfrak{f}': \mathcal{X}'\to \mathcal{Z}$.
Thus we have a bimeromorphic map $\mathfrak{f}: \mathcal{X}\dashrightarrow \mathcal{Z}$.
By Zariski's Main Theorem this map is actually a proper morphism.
This proves Theorem \ref{theorem-main-def}.
\end{proof}

\section{Case \typ{cD/3}}
In this section we prove Theorem \ref{theorem-main-cD/3}.
\begin{de} \textbf{Setup.}
\label{setup-cD/3}
Let $(X,C)$ be an extremal curve germ 
and let $f: (X,C)\to (Z,o)$ be the corresponding contraction.
In particular, $f$ can be flipping.
Throughout this section we assume that $(X,C)$ is of type \type{IA} and 
the only non-Gorenstein point $P\in (X,C)$ is of type \typ{cD/3}
\cite{Mori-1985-cla}, \cite{Reid-YPG1987}.
Our arguments here are very similar to those in 
\cite[\S 6]{Kollar-Mori-1992}.
Note that by Corollary \ref{Corollary-T-singularities}
the point $(H,P)$ is not log terminal for any divisor $H\in |\OOO_X|_C$.
Let $\sigma=(\sigma_1,\dots,\sigma_n)$ be a weight. 
Below, for a formal power series
$\alpha$ in $n$ variables,
$\alpha_{\sigma=m}$ means the sum of the monomials in $\alpha$
whose $\sigma$-weight is $m$.
Put $\sigma:=(1,1,2,3)$.
As in \cite[6.5]{Kollar-Mori-1992}, up to coordinate 
change the point $(X,P)$ is given by
\[
\{\alpha (y_1,y_2,y_3,y_4)=0\}\subset\CC^4_{y_1,y_2,y_3,y_4}/\muu_3(1,1,2,0),
\]
where 
\[
\alpha=y_4^2+y_3^3+\delta_3(y_1,y_2)
+(\text{terms of degree $\ge 4$}),
\]
$\delta_3(y_1,y_2)=\alpha_{\sigma=3}(y_1,y_2,0,0)\neq 0$, 
$\wt \alpha \equiv 0\mod 3$, and $C^\sharp$ is the $y_1$-axis.
If $\delta_3(y_1,y_2)$ is squarefree (resp. has a double factor, is a cube of a linear form),
then $(X,P)$ is said to be a \textit{simple} (resp. \textit{double}, \textit{triple}) \typ{cD/3} point.
The general member $F\in |{-}K_X|$ 
modulo a coordinate change 
is given by the equation $y_1=0$
(see \cite{Reid-YPG1987}). 

\end{de}

\begin{thm} {\bf Lemma.}
\label{lemma-equation-H}
In the above coordinate system 
there exists a member $H\in |\OOO_X|_C$
given by the equation $y_4=\xi $, where $\xi=\xi(y_1,\, y_2,\, y_3)$ is an invariant in
the ideal $(y_2,\, y_3)^3+y_1(y_2,y_3)$.
\end{thm}

\begin{proof}
We have the following exact sequence
\[
0\longrightarrow \omega_X \longrightarrow \OOO_X 
\longrightarrow \OOO_F \longrightarrow 0.
\]
If $f$ is a birational contraction, then 
$R^1f_*\omega_X=0$ by the Grauert-Riemenshneider vanishing theorem.
Hence any section $\bar s\in \OOO_F$
lifts to a section $s\in f_*\OOO_X$.
So, the assertion is clear in this case.
Assume that $f$ is a $\QQ$-conic bundle.
Obviously, $\tau:=f|_F$ is a double cover. 
Since $R^1f_*\omega_X=\omega_Z$ (see \cite[Lemma 4.1]{Mori-Prokhorov-2008})
and $\omega_F\simeq\OOO_F$, we have
\[
f_*\OOO_X \longrightarrow \tau_* \omega_F \longrightarrow
\omega_Z\longrightarrow 0.
\]
The last map is nothing but the trace map $\mt{Tr}_{F/Z}: \tau_*\omega_F\to \omega_Z$.
According to \cite[2.1-2.2]{Mori-Prokhorov-2008III} the induced map 
\[
f_*\OOO_X \longrightarrow \tau_* \omega_F/\tau^*\omega_Z
\]
is surjective. 
We may assume that the equation of $F$
in $\CC^3_{y_2,y_3,y_4}$ is as follows 
\[
\beta(y_2,y_3,y_4):=\alpha(0,y_2,y_3,y_4)=y_4^2+y_3^3+\delta_3(0,y_2)+(\text{terms of degree $\ge 4$}).
\]
Locally, near $P^\sharp$, the sheaf $\omega_{F^\sharp}$ is generated by
\[
\eta:=
\mt{Res}\frac{dy_2\wedge d y_3\wedge d y_4}{\beta}
= 
-\frac{dy_2\wedge d y_3}{\partial \beta/\partial y_4}
=
\frac{dy_2\wedge d y_4}{\partial \beta/\partial y_3}
=
-\frac{dy_3\wedge d y_4}{\partial \beta/\partial y_2}.
\]
Since $\eta$ is an invariant, it is also a generator of $\omega_{F}$
near $P$. 
Further, since $Z$ is smooth, one has
\[
\tau^*\Omega_Z^2=\tau^*\omega_Z\subset \Omega_F^2\longrightarrow
\omega_F.
\]
The generators of $\OOO_{F,P}$ are $y_4$, $w:=y_2y_3$, $u:=y_2^3$,
and $v:=y_3^3$ with relations $uv=w^3$ and $y_4^2+v+u+\cdots=0$. 
Eliminating $v$ we get three generators $y_4$, $w$, $u$
and one relation $u(u+y_4^2+\cdots)+w^3=0$.
Hence 
$\Omega_F^2$ is generated by the elements
\[
\begin{array}{lllll}
d w\wedge d u&=&d(y_2y_3)\wedge d(y_2^3)&= & 3 y_2^3 d y_3\wedge d y_2, 
\\[5pt]
d u\wedge d y_4&=& d(y_2^3)\wedge d y_4&=& 3y_2^2 d y_2\wedge d y_4, 
\\[5pt]
d w\wedge d y_4&=&d(y_2y_3)\wedge d y_4&=& y_2 d y_3\wedge d y_4+y_3 d y_2\wedge d y_4.
\end{array}
\]
Then $\Omega_F^2$ is contained in $\eta I$, where
\[
I:=\langle 
y_3^3\partial \beta/\partial y_4,\
y_2^2\partial \beta/\partial y_3,\
y_2\partial \beta/\partial y_2,\
y_3\partial \beta/\partial y_3
\rangle\subset (y_2,\, y_3,\, y_4)^3.
\]
So, $\tau^*\omega_Z\subset (\tau_*\omega_F)I$.
Therefore, for some $\xi\in I$ the section $\bar s=y_4-\xi \in \OOO_F$
lifts to a section $s\in f_*\OOO_X$.
Since 
\[
s \equiv y_4 \mod (y_2,\, y_3,\, y_4)^3+y_1(y_2,y_3,\, y_4),
\]
one can apply Weierstrass' preparation theorem to get Lemma \ref{lemma-equation-H}.
\end{proof}

\begin{thm} {\bf Corollary.}
\label{corollary-H-normal}
If $y_4$ is a part of an $\ell$-free $\ell$-basis of $\gr_C^1\OOO$, then 
a general member $H\in |\OOO_X|_C$ is normal.
\end{thm}

\begin{de}
\label{par-def-ell}
Recall that $\ell(P):=\len_{P^\sharp} I^{\sharp (2)}/I^{\sharp 2}$
\cite[9.4.7]{Mori-1988}. According to \cite[2.16]{Mori-1988}
we have $i_P(1)=\down{\ell(P)/3}+1$ and the coordinate system
$(y_i)$ can be chosen so that 
$\alpha \equiv y_1^{\ell(P)}y_i \mod (y_2,y_3,y_4)^2$, where 
$i\in \{2,\,3,\,4\}$ and $\ell(P)+\wt y_i\equiv 0\mod 3$.
Since $(X,P)$ is of type \typ{cD/3}, we have $\ell(P)>1$.
\end{de}

Now we are going to prove
Theorem \ref{theorem-main-cD/3}
by considering cases according to the value of $\ell(P)$.
We start with the case $\ell(P)=2$.

\begin{thm}{\bf Theorem.}\label{(6.2)} 
Let the notation and assumptions be as in \xref{setup-cD/3}.
Assume that $\ell(P)=2$ or, equivalently, $i_P(1)=1$. Then the following assertions hold.

\begin{sthm}\label{(6.2.0)}
The contraction $f$ is birational; the general member $H \in |\OOO_X|_C$ 
and its image $T =f(H) \in |\OOO_Z|$ are normal and have only rational
singularities.
\end{sthm}
\begin{sthm}{}
\label{(6.2.3)} 
If $f$ is flipping \textup(resp. divisorial\textup), then 
$P$ is not a triple \typ{cD/3} point and
the dual graph of $(H,C)$ is given as follows with $a=0$ \textup(resp. $a=1$\textup).
\newpar{\label{(6.2.3.1)}
Case of simple \typ{cD/3} point $P$\textup:
\begin{equation*}
\begin{array}{lc@{\,}c@{\,}c}
&
\underbrace{\circ\lin\cdots\lin\circ}_{a}
\lin\bullet\lin\stackrel{3}{\circ}\lin&\circ&\lin\stackrel{3}{\circ}
\\[-15pt]
&&|&
\\[-4pt]
&&\stackunder{3}{\circ}&
\end{array} 
\end{equation*}
}

\newpar{\label{(6.2.3.2)}
Case of double \typ{cD/3} point $P$\textup:
\begin{equation*}
\begin{array}{lc@{\,}c@{\,}c}
&&\stackrel{}{\circ}&
\\[-4pt]
&&|&
\\[-4pt]
&\underbrace{\circ\lin\cdots\lin\circ}_{a}
\lin\bullet\lin\stackunder{3}{\circ}
\lin\stackunder{}{\circ}\lin&\stackunder{3}\circ
&\lin\stackunder{}{\circ}
\\[-7pt]
&&|&
\\[-4pt]
&&\stackrel{}{\circ}&
\end{array} 
\end{equation*}
}
\end{sthm}

\begin{sthm}{}
\label{(6.5.3-6.2.4)}
$\gr^1_C \OOO = (a) \toplus (-a+P^\sharp)$
\end{sthm}
\end{thm}

We now start the proof of Theorem \ref{(6.2)}.

\begin{proof}
Additionally to \xref{setup-cD/3} we assume that $\ell(P)=2$.
Then by \cite[2.16]{Mori-1988} $i_P(1)=1$ and (in some coordinate system) 
$\alpha$ satisfies 
$\alpha\equiv y_1^2y_2\mod (y_2,y_3,y_4)^2$.
Here $C^\sharp$ is the $y_1$-axis as above.
Hence, $y_3$, $y_4$ form an $\ell$-basis of $\gr_C^1\OOO$.
By Corollary \ref{corollary-H-normal} 
$H$ is normal and by Lemma \ref{lemma-surfaces-connect}
$H\setminus \{P\}$ can have 
at most one singular point $R$ which is Du Val.
Therefore,
$X$ can have at most one type \type{III} point.

\begin{sde}
\label{(1.1)} 
\textbf{Subcase $\alpha_{\sigma=3}(y_1,y_2,0,0)$ is squarefree.}
By \cite[6.7.1]{Kollar-Mori-1992} and Lemma \ref{lemma-surfaces-connect} 
the graph $\Delta(H,C)$ is of the form 
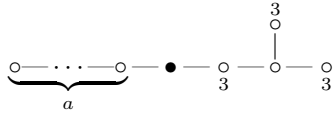
\[
\begin{array}{c@{\,}c@{\,}c}
&\stackrel{3}{\circ}&
\\[-4pt]
&|&
\\[-4pt]
\underbrace{\circ\lin\cdots\lin\circ}_{a}\lin\bullet\lin
\stackunder{3}{\circ}\lin&\circ&\lin\stackunder{3}{\circ}
\end{array} 
\]
We have $a\le 1$, since the corresponding matrix is negative semi-definite. 
But then this matrix is negative definite.
Hence the contraction $f$ is birational. 
If $a=1$, then $H$ is contracted to a singularity $T=f(H)$ of type \typ{A_2}.
Since $T$ is Gorenstein, $f$ is a divisorial contraction
as in \ref{eq-theorem-main-cD/3-A2}.
If $a=0$, that is, $P$ is the only singular point of $H$, 
then $H$ is contracted to a singularity $T=f(H)$
with the following dual graph 
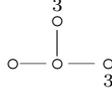
\[
\begin{array}{c@{\,}c@{\,}c}
&\stackrel{3}{\circ}&
\\[-4pt]
&|&
\\[-4pt]
\circ\lin&\circ&\lin\stackunder{3}{\circ}
\end{array} 
\]
Let $s \in H^0(X,\OOO_X)$ be the section defining $H$.
Then $s\OOO_C \subset \gr^1_C \OOO$ is a subbundle
outside $P$ since $H \setminus \{P\}$ is smooth.
At $P^\sharp$, $s\OOO_C^\sharp$ is a subbundle of $\gr^1_C \OOO^\sharp$
by Lemma \ref{lemma-equation-H}, whence $s\OOO_C \simeq (0)$ with $\ell$-structure.
Since $\deg \gr^1_C \OOO = 0$ by $i_P(1)=0$, we have
$\gr^1_C \OOO = (0) \toplus (P^\sharp)$.
Thus $f$ is flipping by \cite[(6.2.4)]{Kollar-Mori-1992}).
\end{sde}

By \ref{(1.1)} it remains to consider the case where 
$\alpha_{\sigma=3}(y_1,y_2,0,0)$
has a double factor.
Note that $y_2$ divides $\alpha_{\sigma=3}(y_1,y_2,0,0)$
because $C^\sharp =(\text{$y_1$-axis})\subset X^\sharp$.
Since $\ell(P)=2$, $y_2y_1^2\in \alpha$.
Then making a coordinate change 
$y_1 \mapsto y_1+cy_2$ we get $\alpha_{\sigma=3}(y_1,y_2,0,0)=y_1^2y_2$
and $C^\sharp$ unchanged.

\begin{sde}
\label{(1.2)} 
\textbf{Subcase $\alpha_{\sigma=3}(y_1,y_2,0,0)=y_1^2y_2$ and
$\alpha_{\sigma=6}(0,y_2,y_3,0)$ is squarefree.}
As above by \cite[6.7.2]{Kollar-Mori-1992} and Lemma \ref{lemma-surfaces-connect} 
the graph $\Delta(H,C)$ is of the form 
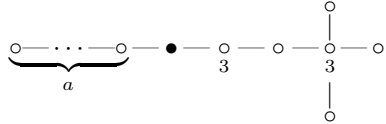
\[
\begin{array}{c@{\,}c@{\,}c@{\,}c}
&\stackrel{}{\circ}&
\\[-4pt]
&|&
\\[-4pt]
\underbrace{\circ\lin\cdots\lin\circ}_{a}\lin\bullet\lin\stackunder{3}{\circ}
\lin\stackunder{}{\circ}\lin&\stackunder{3}\circ
&\lin&\circ
\\[-5pt]
&|&
\\[-4pt]
&\stackrel{}{\circ}&
\end{array} 
\]
with $a\le 1$. 
Again if $a=1$, then $T$ is Du Val of type \typ{D_4}, 
so $f$ is a divisorial contraction as in \ref{eq-theorem-main-cD/3-D4}.
If $a=0$, then similarly to \ref{(1.1)} the contraction $f$ is flipping
(cf \cite[(6.2.3.2)]{Kollar-Mori-1992}).
Since $s\OOO_{C}^{\sharp}$ is a subbundle of 
$\gr_{C^\sharp}^1\OOO_C^\sharp$ at $P^\sharp$,
as we saw above, it is easy to see \ref{(6.5.3-6.2.4)}.
\end{sde}

\begin{sde}
\label{(1.3)-(1.4)} 
\textbf{Subcase $\alpha_{\sigma=3}(y_1,y_2,0,0)=y_1^2y_2$, and
$\alpha_{\sigma=6}(0,y_2,y_3,0)$
has a multiple factor.}
We will show that this case does not occur.
Assume that $f$ is birational. 
Then the map $H^0(\OOO_X)\to H^0(\OOO_F)$ is surjective
\cite[Th. 1.2]{Mori-1988}. Therefore, for any $\lambda\in \CC^*$ there is 
a semi-invariant $\delta$ with $\wt \delta=2$ such that the section 
$y_4+\lambda y_2^3 + \delta y_1$ extends to some element $H'\in |\OOO_X|_C$. 
After the coordinate change $y_4'= y_4+\lambda y_2^3 + \delta y_1$
we see that $H'$ is given by $y_4'=0$ and
$\alpha'=\alpha(y_1,y_2,y_3,y_4'-\lambda y_2^3 - \delta y_1)$.
Note that $y_4^2\in \alpha$, $y_4\notin \alpha$, and 
$\alpha$ may contain $y_2^3y_4$.
Thus $\alpha'_{\sigma=3}(y_1,y_2,y_3,0)= \alpha_{\sigma=3}(y_1,y_2,y_3,0)$
and $\alpha'_{\sigma=6}(0,y_2,y_3,0)= \alpha_{\sigma=6}(0,y_2,y_3,0)+(\lambda^2+c\lambda)y_2^6$
for some $c\in \CC$.
Hence we may assume that $\alpha_{\sigma=6}(0,y_2,y_3,0)$ is square free.
This contradicts our assumption.
(In fact, the above arguments show that the chosen $H$ is not general).

Therefore, $f$ is a $\QQ$-conic bundle.
By Lemma \ref{lemma-H} $(H,P)$ is a rational singularity and 
by Lemma \ref{lemma-equation-H} this 
singularity is analytically isomorphic to 
\[
\{\gamma(y_1,y_2,y_3)=0\}/\muu_3(1,1,2),
\]
where 
$\gamma(y_1,y_2,y_3):=\alpha(y_1,y_2,y_3,\xi)$,
and $C\subset H$ is the image of $y_1$-axis.
Note that the pair $(H,C)$ is not plt at $P$. Indeed, 
otherwise the singularity $\{\gamma=0\}$ is log terminal
(see \cite[Cor. 20.4]{Utah}). Hence it is Du Val.
On the other hand, $\ord \gamma>2$, a contradiction.
Let $\sigma':=(1,1,2)$.
Note that $\gamma_{\sigma'=6}(0,0,1)\neq 0$ because 
$y_3^3\in \alpha$.
Consider the weighted $\sigma'$-blowup 
$\varsigma: \ov H\subset \ov {\CC^3/\muu_3}
\to H\subset \CC^3/\muu_3$. 
Let $\Xi:=\varsigma^{-1}(0)_{\red}$.
The exceptional divisor $\Theta\subset \ov H$ is given in 
$\Xi\simeq \PP(1,1,2)$ by the equation 
$\gamma_{\sigma'=3}(y_1,y_2,0)=y_1^2y_2=0$.
Hence, $\Theta=2\Theta_1+\Theta_2$, where $\Theta_i$ are 
irreducible toric divisors in $\PP(1,1,2)$.
The proper transform $\ov C$ of $C$ meets $\Xi\simeq \PP(1,1,2)$
at the point $\{y_2=y_3=0\}$. So, $\ov C\cap \Theta_1=\emptyset$.
Since $\Theta_2$ is a smooth reduced component of the Cartier divisor
$\Theta=\Xi\cap \ov H$ on $\ov H$, we see that
$\ov H$ is smooth at points on $\Theta_2\setminus \Theta_1$. 

In the chart $U_3\simeq \CC^3/\muu_2(1,1,1)$ over $\{y_3\neq 0\}$ we have a new coordinate system
$y_1 \mapsto y_1y_3^{1/3}$, $y_2 \mapsto y_2y_3^{1/3}$, $y_3 \mapsto y_3^{2/3}$.
Here the surface $\ov H$ is given by the equation 
$y_1^2y_2+\gamma_{\sigma'=6}(y_1,y_2,1)y_3+(\cdots)y_3^2=0$, where 
$\gamma_{\sigma'=6}(0,0,1)\neq 0$. The origin $O_3\in \ov H\cap U_3$ is a 
Du Val point of type \typ{A_1}.
Components $\Theta_1$ and $\Theta_2$ of the exceptional divisor
meet each other at $O_3$ and the pair $(\ov H,\Theta_1+\Theta_2)$ is LC at $O_3$. 
Outside of $O_3$, \ $\ov H$ is a hypersurface and has only rational singularities.
Therefore, the singularities of $\ov H$ are Du Val.
Thus the configuration of curves $\ov C$, $\Theta_1$, and $\Theta_2$ on $\ov H$ 
looks as follows: 
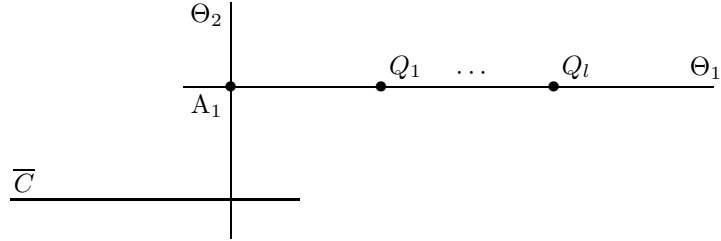
\[
\xy
(-20,0)="A" *{}*+!DR{\ov C},
(10,0)="B" *{}*+!DR{},
(5,-1)="C" *{}*+!DR{},
(5,22)="D" *{}*+!DR{\Theta_2},
(3,15)="E" *{}*+!DL{},
(65,15)="F" *{}*+!DL{\Theta_1},
(25,15)="EE" *{\bullet}*+!DL{Q_1},
(37,15)="EE" *{}*+!D{\cdots},
(25,15)="TT" *{}*+!DL{},
(48,15)="FF" *{\bullet}*+!DL{Q_l},
{"E";"F":"C";"D",x} ="I" *{\bullet}*+!UR{\mathrm A_1},
"B";"A"**{} +/1pc/;-/1pc/ **@{-},
"C";"D"**{} +/1pc/;-/1pc/ **@{-},
"E";"F"**{} +/1pc/;-/1pc/ **@{-}
\endxy
\]
where $Q_1,\dots,Q_l$ are some Du Val points.
By Lemma \ref{lemma-surfaces-connect}
the dual graph $\Delta(H,C)$ is of the form
\begin{newequation}
\label{eqnumer-1}
\begin{array}{c@{\,}c@{\,}c@{\,}c@{\,}c@{\,}c@{\,}}
\overbrace{\circ\lin \cdots\lin\circ}^a\lin
\stackunder{C}{\bullet}
\lin
\stackud{\Theta_2}{b_2}{\circ}
\lin \circ\lin&
\stackud{\Theta_1}{b_1}{\circ}
&\lin&\vdots\lin\whitebox
\end{array}
\end{newequation}
where the box on the right hand side 
indicates some Du Val graphs and the number 
of these Du Val tails is not important.
This configuration forms a fiber of a rational curve fibration.
Contracting black vertices successively we obtain 
\begin{newequation}
 \label{eqnumer-2}
\begin{array}{c@{\,}c@{\,}c@{\,}c@{\,}c@{\,}}
\stackud{\Theta_2}{b_2-a-1}{\circ}
\lin \circ\lin&
\stackud{\Theta_1}{b_1}{\circ}\lin
&\vdots\lin\whitebox
\end{array}
\end{newequation}
This is again a dual graph of a fiber of a rational curve fibration.
Hence $b_2-a-1=1$ and we further obtain 
\[
\begin{array}{c@{\,}c@{\,}c@{\,}c@{\,}c@{\,}}
&
{\scriptstyle{b_1-1}}\ \vspace{5pt} {\circ}
&\lin
&\vdots\lin\whitebox
\end{array}
\]
Hence $b_1=2$ (because the last graph must contain a $(-1)$-vertex) and so the graph 
\ref{eqnumer-2} consists of $(-2)$ and $(-1)$-curves.
Furthermore the graph \ref{eqnumer-2} is not a linear chain because the pair
$(H,C)$ is not PLT at $P$.
In this situation there is only one possibility (see, e.g., \cite[Lemma 7.1.3, 7.1.12]{Prokhorov-2001}):
\[
\begin{array}{c@{\,}c@{\,}c@{\,}c@{\,}c@{\,}c@{\,}c@{\,}c@{\,}}
&&&&&&&\circ
\\[-5pt]
&&&&&&\diagup
\\[-5pt]
\stackunder{\Theta_2}{\bullet}
\lin \circ\lin&
\stackunder{\Theta_1}{\circ}
&\lin&\cdots&\lin&\circ
\\[-8pt]
&&&&&&\diagdown
\\[-5pt]
&&&&&&&\circ
\end{array}
\]
Therefore, the original graph \ref{eqnumer-1} is of the form
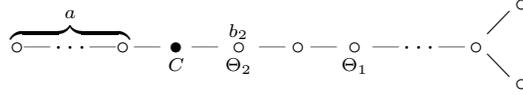
\[
\begin{array}{c@{\,}c@{\,}c@{\,}c@{\,}c@{\,}c@{\,}c@{\,}c@{\,}}
&&&&&&&\circ
\\[-5pt]
&&&&&&\diagup
\\[-12pt]
\overbrace{\circ\lin \cdots\lin\circ}^a\lin
\stackunder{C}{\bullet}
\lin
\stackud{\Theta_2}{b_2}{\circ}
\lin \circ\lin&
\stackunder{\Theta_1}{\circ}
&\lin&\cdots&\lin&\circ
\\[-10pt]
&&&&&&\diagdown
\\[-5pt]
&&&&&&&\circ
\end{array}
\]
But then $H$ has only log terminal singularities
(see, e.g., \cite[Ch. 3]{Utah}).
Hence $H$ has only \typ{T}-singularities (see \ref{T-singularities-definition}) while the 
right hand side singularity is not of type \typ{T}
(see Proposition \ref{T-singularities-characterization}), a contradiction.
Thus the case \ref{(1.3)-(1.4)} does not occur. 
\end{sde}
Now the assertion of Theorem \ref{(6.2)} follows from 
\ref{(1.1)}, \ref{(1.2)}, and \ref{(1.3)-(1.4)}. 
This completes our treatment of the case $\ell(P)=2$.
\end{proof}

\begin{thm} {\bf Corollary.}
\label{corollary-3-III}
In the notation of \xref{setup-cD/3}
$X$ has at most one type \type{III} point.
\end{thm}

\begin{proof}
If $X$ has two type \type{III} points $R_1$ and $R_2$, then
by \cite[(2.3.3)]{Mori-1988} and \cite[(3.1.5)]{Mori-Prokhorov-2008}
we have $i_P(1)=i_{R_1}(1)=i_{R_2}(1)=1$. Then by \cite[2.16]{Mori-1988}
$\ell(P)=2$. This contradicts Theorem \ref{(6.2)}. 
\end{proof}

\begin{thm} {\bf Lemma (cf. {\cite[Lemma 6.12]{Kollar-Mori-1992}}).}
\label{lemma-bound-ell}
If, in the notation of \xref{setup-cD/3},
$X$ has a type \type{III} point, then $\ell(P)\le 4$ and 
$i_P(1)\le 2$. 
\end{thm}

\begin{proof}
Assume that $\ell(P)\ge 5$.
As in \ref{par-def-ell} take a coordinate system 
so that 
$\alpha \equiv y_1^{\ell(P)}y_i \mod (y_2,y_3,y_4)^2$, where 
$i\in \{2,\,3,\,4\}$ and $\ell(P)+\wt y_i\equiv 0\mod 3$.
Similarly to the proof of \cite[Lemma 6.12]{Kollar-Mori-1992}
we use the deformation $\alpha_\lambda=\alpha+\lambda y_1^{\ell(P)-3}y_i$
(see Theorem \ref{theorem-main-def})
and get a germ $(X_{\lambda},C_{\lambda})$ with two type \type{III}
points and a point of type \typ{cD/3}. 
This contradicts Corollary \ref{corollary-3-III}.
\end{proof}

For the case $\ell(P) \ge 3$, we are going to prove the following,
which settles Theorem \ref{theorem-main-cD/3}.

\begin{thm}{\bf Theorem.}
\label{(6.3)}
Let the notation and assumptions be as in \xref{setup-cD/3}.
Assume $\ell(P) \ge 3$ or, equivalently, $i_P(1)\ge 2$. 
Then the following assertions hold.

\begin{sthm}{}
\label{(6.3.a)}
$\ell(P)=3$ or $4$ \textup(i.e. $i_P(1)= 2$\textup), and $f$ is birational.
\end{sthm}
\begin{sthm}{}\label{(6.3.0)}
$P$ is a double \textup(resp. triple\textup) \typ{cD/3} point if $(X,C)$ 
is isolated \textup(resp. divisorial\textup).
\end{sthm}

\begin{sthm}{}\label{(6.3.1)}
$ X $ is smooth outside of $P$ and there is an $\ell $-isomorphism 
\begin{newequation}
\label{(6.3.1.1)}
\gr^{1}_{C}\OOO =
((4-\ell(P))P^\sharp) \toplus (-1+2P^\sharp).
\end{newequation} 
\end{sthm}

\begin{sthm}{}\label{cD3-ip=2-C=ci}
For general members $D\in |K_X|$ and $D' \in |K_X|$ 
\textup(resp. $D' \in |\OOO_X|_C$\textup),
$D \cap D'$ is equal to $4C$ \textup(resp. $3C$\textup) as a $1$-cycle.
\end{sthm}

\begin{sthm}{}\label{(6.3.3)} 
The general member $H \in |\OOO_X|_C$ and
its image $T =f(H) \in |\OOO_Z|$ are normal and have only rational
singularities. 
The dual configuration of $(H,C)$ is as follows. 

\newpar{
Case of isolated $(X,C)$\textup:
\begin{equation*}\label{(6.2.3.1-a)}
\begin{array}{lc@{\,}c@{\,}l}
&&\circ&\\
[-4pt]
&&|\\
&\bullet\lin\circ
\lin&\stackunder{3}{\circ}&\lin\circ\lin\stackunder{3}\circ
\\
&&|
\\[-4pt]
&&\circ
\end{array} 
\end{equation*}
}

\newpar{
Case of divisorial $(X,C)$\textup:
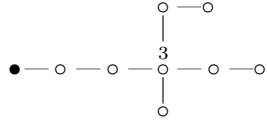
\begin{equation*}\label{(6.2.3.1-b)}
\begin{array}{lc@{\,}c@{\,}l}
&&\circ&\lin\circ\\
[-4pt]
&&|\\
&\bullet\lin\circ\lin\circ
\lin&\stackrel{3}{\circ}&\lin\circ\lin\circ\\
[-4pt]
&&|
\\[-4pt]
&&\circ
\end{array}
\end{equation*}
}
\end{sthm}

\begin{sthm}{}\label{converse-cD3-ip=2}
Conversely if $(X,C)$ is an arbitrary germ of a threefold along $C \simeq \PP^1$
with a double \textup(resp. triple\textup) \typ{cD/3} point $P\in C$.
If $(X,C)$ satisfies \xref{(6.3.1)}, then $(X,C)$ is an isolated \textup(resp. a divisorial\textup)
extremal curve germ.
\end{sthm}
\end{thm}

\begin{proof}
In the hypothesis of \xref{setup-cD/3} we additionally assume that 
$\ell(P) \ge 3$.

\begin{thm}{\bf Lemma.}
\label{(2)} 
Under notation of \xref{(6.3)} $X$ has no type \type{III} points.
\end{thm}

\begin{proof}
Assume that $X$ has a type \type{III} point $R$. We derive a contradiction.
By Lemma \ref{lemma-bound-ell}\quad 
$\ell(P) =3$ or $4$.
\begin{sde}
\label{(2)-ell=3} 
\textbf{Case $\ell(P)=3$.}
We claim that $H^1(\gr_C^2\omega)\neq 0$.
By \cite[2.16]{Mori-1988}\ $i_P(1)=2$ and (in some coordinate system) 
$\alpha$ satisfies $\alpha\equiv y_1^3y_4\mod (y_2,y_3,y_4)^2$
(and $C^\sharp$ is the $y_1$-axis). If $\alpha$ contains the term
$y_1^ky_2y_3$, then $k\ge 3$ and this term can be removed by 
the coordinate change $y_4\mapsto y_4-y_1^{k-3}y_2y_3$.
Hence, we may assume that
\begin{equation*}
\alpha\equiv y_1^3y_4+
\lambda y_1y_2^2+\mu y_1^2y_3^2
\mod (y_2,y_3)^3+y_4(y_2,y_3,y_4)\quad \bigl(\subset I_C^{(3)\sharp}\bigr)
\end{equation*}
for some $\lambda,\, \mu\in \OOO_X \mod I_C$.
The functions $y_2$, $y_3$ form an $\ell$-basis of $\gr_C^1\OOO$ at $P$.
Since 
\[
\deg \gr_C^1\OOO=1-i_P(1)-i_R(1)=-2
\] 
and $H^1(\gr_C^1\OOO)=0$, we have
$\gr_C^1\OOO=\OOO(-1)\oplus \OOO(-1)$.
Furthermore, 
by \cite[(2.8)]{Kollar-Mori-1992}, one has
\[
\gr_C^1\OOO=
 (-1+P^\sharp)\toplus\ (-1+2P^\sharp).
\]
where $y_3$ (resp. $y_2$) is an $\ell$-free $\ell$-basis of $(-1+P^\sharp)$
\textup(resp. of $(-1+2P^\sharp)$\textup) at $P$.
Let $\sigma$ be an $\ell$-basis of $\omega$. 
By the above, one has
\[
\omega\totimes \ti S^2 \gr_C^1\OOO =
(-2+P^\sharp)\toplus (-2+2P^\sharp) \toplus (-1),
\]
where $y_3^2\sigma$ \textup(resp. $y_3y_2\sigma$, $y_2^2\sigma$\textup) 
is an $\ell$-free $\ell$-basis of $(-2+P^\sharp)$
\textup(resp. $(-2+2P^\sharp)$, $(-1)$\textup) at $P$.
There is an injection of coherent sheaves
\[
\iota: 
\omega\totimes \ti S^2 \gr_C^1\OOO \longrightarrow \gr_C^2\omega.
\]
As an abstract sheaf, $\omega\totimes \ti S^2 \gr_C^1\OOO$ at $P$
is generated by sections $y_3^2y_1\sigma$, $y_3y_2y_1^2\sigma$, $y_2^2\sigma$.
Further, it is easy to see that 
$I_C^{(2)\sharp}/I_C^{(3)\sharp}$ at $P$ is generated by 
elements $y_4$, $y_3^2$, $y_2y_3$, $y_2^2$. Hence,
$\gr_C^2\omega$ at $P$ is generated by 
$y_4y_1^2\sigma$, $y_3^2y_1\sigma$, $y_2y_3y_1^2\sigma$, $y_2^2\sigma$.
On the other hand, $y_4\in I_C^{(2)}$ and $y_1^2y_2y_4,\, y_1y_3y_4\in I_C^{(3)}$.
By our expression of $\alpha$ 
\[
(y_1^2y_4+\lambda y_2^2+\mu y_1y_3^2)\sigma =0 \text{\quad in $\gr_C^2\omega$ at $P$.}
\]
Hence, 
$\gr_C^2\omega$ at $P$ is generated by the elements
$y_3^2y_1\sigma$, $y_2y_3y_1^2\sigma$, $y_2^2\sigma$.
This means that $\iota$ is an isomorphism at $P$.

Since $i_R(1)=1$, by \cite[2.16]{Mori-1988} $\ell(R)=1$ and in some coordinate 
system the local equation 
$\beta(z_1,\dots,z_4)=0$ of $(X,R)$ satisfies 
$\beta \equiv z_1z_2\mod (z_2,z_3,z_4)^2$, where $C$ is the $z_1$-axis.
Then locally near $R$ we have $I_C^{(2)}=(z_3^2, z_3z_4, z_4^2, z_2)$, so 
\[
\gr_C^1\OOO=\OOO z_3\oplus \OOO z_4\quad \text{and} \quad
S^2\gr_C^1\OOO= \OOO z_3^2\oplus \OOO z_4^2\oplus \OOO z_3z_4.
\]
Furthermore, $\gr_C^2\OOO$ is generated by 
$z_2$, $z_3^2$, $z_4^2$, $z_3z_4$. Hence, $z_2$ generates
$\Coker \iota$ and so $\len_R \Coker \iota\le 1$.
In this case,
$\dim H^0(\Coker \iota)\le 1$ and $\dim H^1(\im \iota)=2$.
Therefore, $H^1(\gr_C^2\omega)\neq 0$ as claimed. 

Now from 
$H^0(\gr_C^j\omega)=0$, where $j=0,\, 1$ and 
the exact sequences
\[
0\longrightarrow \gr_C^n \omega \longrightarrow \omega_X/F^{n+1}\omega_X
\longrightarrow \omega_X/F^n\omega_X \longrightarrow 0, \quad n=1,\, 2
\]
we have $H^1(\omega_X/F^3\omega_X)\neq 0$. 
If $f$ is birational, then by \cite[1.2-1.2.1]{Mori-1988} we get a contradiction.
Assume that $f$ is a $\QQ$-conic bundle.
Put $V:=\Spec_X \OOO_X/ I_C^{(3)}$.
By \cite[Th. 4.4]{Mori-Prokhorov-2008} $V\supset f^{-1}(o)$.
Since 
\[
-K_X\cdot V=-6K_X\cdot C=2=-K_X\cdot f^{-1}(o), 
\]
we have $V= f^{-1}(o)$.
Let $P\in C$ be a
general point. Then in a suitable coordinate system $(x,y,z)$ near
$P$ we may assume that $C$ is the $z$-axis. So, $I_C=(x,y)$ and
$I_C^{(3)}=(x^3,x^2y, xy^2, y^3)$. But then $V= f^{-1}(o)$ is not a local complete
intersection near $P$, a contradiction.
This disproves the case \ref{(2)-ell=3}.
\end{sde}

\begin{sde}
\label{(2)-ell=4} 
\textbf{Case $\ell(P)=4$.}
By deformation $\alpha_\lambda=\alpha+\lambda y_1^3y_4$ at $(X,P)$
we get a germ $(X_\lambda,C_{\lambda})$ with a point 
$P_\lambda$ of type \typ{cD/3} with $\ell(P_\lambda)=3$
(see Theorem \ref{theorem-main-def}).
Moreover, $X_\lambda$ has a point $R_\lambda$ of type 
\type{III}. This is impossible by \ref{(2)-ell=3}.
\end{sde}
This proves \ref{(2)}.
\end{proof}

>From now on we treat the case where $P$ is the only singular point of $X$
and $\ell(P)\ge 3$.

\begin{thm} {\bf Lemma (cf. {\cite[Lemma 6.12]{Kollar-Mori-1992}}).}
\label{lemma-cd3-ip4}
In the notation of \xref{setup-cD/3} we have
$\ell(P)\le 4$ and $i_P(1)\le 2$. 
\end{thm}

\begin{proof}
Assume that $\ell(P)\ge 5$.
Similarly to \cite[Lemma 6.12]{Kollar-Mori-1992} and Lemma \ref{lemma-bound-ell}
we write 
$\alpha \equiv y_1^{\ell(P)}y_j \mod (y_2,y_3,y_4)^2$, where 
$j\in \{2,\,3,\,4\}$ and $\ell(P)+\wt y_j\equiv 0\mod 3$, and use
deformation $\alpha_\lambda=\alpha+\lambda y_1^{\ell(P)-3}y_j$
(see Theorem \ref{theorem-main-def}).
We get a germ $(X_{\lambda},C_{\lambda})$ with a type \type{III}
point $R_\lambda$ and a point $P_\lambda$ of type \typ{cD/3}
with $\ell(P_\lambda)=\ell(P)-3$.
If $\ell(P)\ge 6$, we get a contradiction by the case 
\ref{(2)} considered above.

Hence $\ell(P)=5$, and $X \setminus \{P\}$ is smooth by Lemma \ref{lemma-bound-ell}.
Then $\alpha \equiv y_1^5y_2 \mod (y_2,y_3,y_4)^2$,
$\deg \gr^1_C \OOO_X =-1$, and $y_4, y_3$ form an $\ell$-basis
for $\gr^1_C \OOO_X$.
Thus $H$ is normal at $P$ by Corollary \ref{corollary-H-normal},
and we see that $\gr^1_C \OOO_X = (0) \toplus (-1+P^\sharp)$,
$H$ is smooth outside $P$,
$y_3$ is an $\ell$-basis of $\gr^1_C \OOO_H$,
and $\gr^1_C \OOO_H=(-1+P^\sharp)$.
We also see 
$$
\gr^0_C \omega_H=\gr^0_C \omega_X=(-1+2P^\sharp)
\quad \text{and} \quad\gr^1_C \omega_H=(-1).
$$

We note that $C^\sharp=y_1$-axis 
$\subset H^\sharp \subset \CC_{y_1,y_2,y_3}^3$, 
and $H^\sharp=\{\beta=0\}$,
where 
\[
\beta \equiv y_1^5y_2+c y_1^2y_3^2 
\mod (y_2^2,y_2y_3,y_3^3),
\] 
and $c\in \CC$.
We claim $c\neq0$. Indeed otherwise we have
$y_2 \in \OOO_H(-3C)^\sharp$,
whence $\gr^2_C \OOO_H^\sharp=\OOO_{C^\sharp}y_3^2$ and
$\gr^2_C \OOO_H=(\gr^1_C\OOO_H)^{\totimes 2}=(-2+2P^\sharp)$.
Thus $H^1(H,\OOO_H)\neq0$, a contradiction. Hence $c\neq0$.

Since $P$ is a \typ{cD/3} point, we have $y_2y_3 \not\in \alpha$ and 
$y_3^3 \in \alpha$, and hence $y_2y_3 \not\in \beta$ and $y_3^3 \in \beta$. 
Since $c\neq0$, the terms $\gamma(y_1)y_1^3y_2y_3$ can be
killed by a $\muu_3$-coordinate change $y_3 \mapsto y_3-\gamma(y_1)y_1y_2/(2c)$
and we may further assume
\begin{equation}
\label{3.10.1}
\beta \equiv y_1^5y_2+c y_1^2y_3^2 +y_3^3 \mod (y_2^2,y_2y_3^2,y_3^4).
\end{equation}

We claim that
$\gr^2_C \OOO_H =(-1+2P^\sharp)$ and $\gr^3_C \OOO_H = (-1)$.
First by $y_2 \in \OOO_H(-2C)^\sharp$, one has
$y_1^2(y_1^3y_2+cy_3^2) \in \OOO_H(-3C)^\sharp$.
Hence if we set $z := y_1^3y_2+cy_3^2$,
then $z \in \OOO_H(-3C)^\sharp$ and 
$y_3^2\equiv -y_1^3y_2/c \mod (z)$.
Thus by $\OOO_H(-2C)^\sharp=(y_2,y_3^2)$, we see
$$
\OOO_H(-2C)^\sharp/(y_2^2,y_2y_3,z) = \OOO_{C^\sharp}y_2 \simeq \OOO_{C^\sharp}
\quad\text{and}\quad
\OOO_H(-3C)^\sharp=(y_2^2,y_2y_3,z).
$$
We also have $y_1^2z+y_3^3 \in (y_2^2,y_2y_3^2,y_3^4)$ by (\ref{3.10.1}), 
whence $z\equiv y_1y_2y_3/c \mod (y_2^2,y_2y_3^2,zy_3)$.
Thus
$$
\OOO_H(-3C^\sharp)/(y_2^2,y_2y_3^2,zy_3) = \OOO_{C^\sharp}y_2y_3 \simeq
\OOO_{C^\sharp}
\hspace{5.1pt}\text{and}\hspace{5.2pt}
\OOO_H(-4C)^\sharp=(y_2^2,y_2y_3^2,zy_3).
$$
>From these follows the claim:
$$
\gr^2_C \OOO_H = (\gr^1_C \OOO_H)^{\totimes 2}
(3P^\sharp)=(-1+2P^\sharp)
$$
and 
$$
\gr^3_C \OOO_H = \gr^1_C \OOO_H
\totimes \gr^2_C \OOO_H = (-1).
$$

We then claim $H^1(\omega_H/\omega_H(-4C))\neq0$.
Indeed this follows from
$$
\gr^2_C \omega_H = \gr^0_C\omega_H \totimes \gr^2_C \OOO_H= (-1+P^\sharp)
$$
and 
$$
\gr^3_C \omega_H = \gr^0_C\omega_H \totimes \gr^3_C \OOO_H= (-2+2P^\sharp).
$$

Since $\omega_H=\omega_X \totimes \OOO_H$, the
non-vanishing $H^1(\omega_H/\omega_H(-4C))\neq0$ means that
$f$ is a $\QQ$-conic bundle
\cite[Remark 1.2.1]{Mori-1988}
and the subscheme $4C$ of $H$ contains 
the 
scheme-theoretic fiber $f^{-1}(o)$ \cite[Th. 4.4]{Mori-Prokhorov-2008}. 
However 
$$
-K_X \cdot 4C=4/3<2=-K_X \cdot f^{-1}(o), 
$$
a contradiction. 
The case $\ell(P)=5$ is thus disproved.
\end{proof}

\begin{de}
\label{(4)} 
\textbf{Case $\ell(P)=3$ and no type \type{III} points.}
By \cite[2.16]{Mori-1988} $i_P(1)=2$ and (in some coordinate system) 
$\alpha$ satisfies $\alpha\equiv y_1^3y_4\mod (y_2,y_3,y_4)^2$
(and $C^\sharp$ is the $y_1$-axis).
Hence, $y_2$, $y_3$ form an $\ell$-basis of $\gr_C^1\OOO$.
Since $\deg \gr_C^1\OOO=1-i_P(1)=-1$ and $H^1(\gr_C^1\OOO)=0$, 
$\gr_C^1\OOO=\OOO\oplus \OOO(-1)$.
Further, by \cite[(2.8)]{Kollar-Mori-1992} there are only two possibilities:
\[
\gr_C^1\OOO=
\begin{cases}
\quad (2P^\sharp)&\toplus\ (-1+P^\sharp),
\\
\quad (P^\sharp)&\toplus\ (-1+2P^\sharp).
\end{cases}
\]
Consider the first case, i.e., 
$\gr_C^1\OOO=(2P^\sharp)\toplus (-1+P^\sharp)$.
Then the arguments in the first part of proof of
\cite[(6.13)]{Kollar-Mori-1992} can be applied.
Let $J$ be the $C$-laminal ideal of width $2$ such that
$J/F^2_C\OOO=(2P^\sharp)$. 
Then we conclude that $H^1(\omega/F^4(\omega,J))\neq 0$ 
\cite[p. 599-600]{Kollar-Mori-1992}.
If the contraction $f$ is birational, we get a contradiction by 
\cite[1.2-1.2.1]{Mori-1988}. Let $f$ be a $\QQ$-conic bundle.
Put $V:=\Spec_X \OOO_X/F^4 (\OOO,J)$. Then $V\equiv mC$ for some $m$.
By \cite[Th. 4.4]{Mori-Prokhorov-2008} $V\supset f^{-1}(o)$.
Hence, $m/3=-K_X\cdot V<2=-K_X\cdot f^{-1}(o)$. On the other hand,
near a general point $S\in C$, $J$ is generated by $(z_2, z_3^2)$,
where $(z_1,z_2,z_3)$ are some local coordinates such that $C$ is the
$z_1$-axis. Hence, $F^4(\omega,J)=J^2=(z_2, z_3^2)^2$ near $S$.
So, $m=\len(\CC[z_2,z_3]/F^4(\omega,J))= 6$.
Therefore, $f^{-1}(o)=V$ and its ideal sheaf coincides with $F^4(\omega,J)$. 
However, $F^4(\omega,J)$ is not generated by two elements near $S$,
so $f^{-1}(o)$ is not a locally complete intersection,
a contradiction.

Consider the second case, i.e., $\gr_C^1\OOO=(P^\sharp)\toplus (-1+2P^\sharp)$.
If $(X,P)$ is a double \typ{cD/3} point, then
$f$ is a flipping contraction by 
\cite[Theorem (6.3)]{Kollar-Mori-1992},
whence we get the configuration \ref{(6.2.3.1-a)}.
Thus we assume that the term $y_1y_2^2$ does not appear in $\alpha$.
Further we use arguments of the proof of 
\cite[Lemma (6.13), p. 600]{Kollar-Mori-1992}. 
Let $J$ be the $C$-laminal ideal of width $2$ such that
$J/F^2_C\OOO=(P^\sharp)$. 
Modulo a $\muu_3$-equivariant change of coordinates 
we may further assume that 
$y_3$ (resp. $y_2$)
is an $\ell$-free $\ell$-basis of $(P^\sharp)$ (resp. $(-1+2P^\sharp)$)
in $\gr_C^1\OOO$ and that
$\alpha\equiv y_1^3y_4\mod I^\sharp J^\sharp$.
Whence 
$J^\sharp=(y_2^2, y_3, y_4)$ at $P^\sharp$ and
$y_4\in F^3(\OOO,J)^\sharp$. 
Let $K$ be the ideal such that 
$J\supset K\supset F^3(\OOO,J)$ and $K/F^3(\OOO,J)=(P^\sharp)$ 
in 
\[
\gr^2(\OOO,J)=\gr^{2,0}(\OOO,J)\toplus \gr^{2,1}(\OOO,J)
= (P^\sharp)\toplus (-1+P^\sharp).
\]
Here we may assume that 
$y_3$ (resp. $y_2^2$)
is an $\ell$-free $\ell$-basis of $(P^\sharp)$ (resp. $(-1+P^\sharp)$)
in the above $\ell$-splitting
modulo a coordinate change $y_3 \mapsto y_3+ (\cdots)y_2^2$.
We then have $K^\sharp=(y_2^3,y_3,y_4)$ at $P^\sharp$ and
\[
\gr^1(\OOO,K)=(-1+2P^\sharp),
\qquad
\gr^2(\OOO,K)=(-1+P^\sharp).
\]
We have
$\gr^{3,0}(\OOO,K) \simeq\gr^{2,0}(\OOO,J) \simeq (P^\sharp)$
and
\[
\alpha \equiv y_1^3y_4 + c y_2^3 \mod I^\sharp K^\sharp,
\]
for some unit $c \in \OOO_X^{\times}$, because
$I^\sharp J^\sharp = I^\sharp K^\sharp+(y_2^3)$ and $y_2^3 \in \alpha$. 
Whence we have an
$\ell$-isomorphism
\[
\gr^{3,1}(\OOO,K) \simeq \gr^1(\OOO,K)^{\totimes 3}(3P^\sharp) \simeq (0)
\]
as in \cite[p. 600]{Kollar-Mori-1992}, and an $\ell$-splitting
\[
\gr^3(\OOO,K)=
\gr^{3,0}(\OOO,K)\toplus \gr^{3,1}(\OOO,K),
\]
in which $y_3$ (resp. $y_4$)
is an $\ell$-free $\ell$-basis of $(P^\sharp)$ (resp. $(0)$)
modulo a coordinate change $y_3 \mapsto y_3+ (\cdots) y_1^2y_4$.
For any $l>0$ there is 
a natural exact sequence
\begin{equation}\label{Esq-FlOOOK}
0\longrightarrow F^{l+1}(\OOO,K) \longrightarrow F^{l}(\OOO,K)
\longrightarrow \gr^{l}(\OOO,K) \longrightarrow 0.
\end{equation}
We claim that the sections $y_1y_3,\, y_4\in \gr^3(\OOO,K)$ can be 
extended to sections of $F^3(\OOO,K)= F^1(K)$.
By \eqref{Esq-FlOOOK} it is sufficient to show that $H^1(F^{4}(\OOO,K))=0$.
There are injections of coherent sheaves
\[
\begin{array}{l@{\quad\hookleftarrow\quad}l}
\gr^{3n}(\OOO,K)&\ti S^n\gr^3(K),
\\
\gr^{3n+1}(\OOO,K)&\ti S^n\gr^3(K)
\totimes \gr^1(\OOO,K),
\\
\gr^{3n+2}(\OOO,K)&\ti S^n\gr^3(K)
\totimes \gr^2(\OOO,K)
\end{array}
\]
with cokernels of finite length.
Therefore, for any $l>0$, the degree of each component 
in a decomposition of $\gr^{l}(\OOO,K)$ in a direct sum
is at least $-1$. Then $H^1(\gr^{l}(\OOO,K))=0$ and from 
\eqref{Esq-FlOOOK}
we get surjections 
\[
H^{1}(F^{l+n}(\OOO,K))\to H^1(F^{l}(\OOO,K))\quad
\text{for $l,\, n>0$}.
\] 
Hence, $H^1(F^l(\OOO,K)/F^{l+n}(\OOO,K))=0$.
Note that for any $m>0$ there is $n>0$ such that
$I_C^mF^{l}(\OOO,K)\supset F^{l+n}(\OOO,K)$.
By the Formal Function Theorem we have
\begin{multline*}
H^1(F^l(\OOO,K))\compl=H^1(F^l(\ha{\OOO},K))
=\varprojlim H^1(F^l(\OOO,K)/I_C^mF^{l}(\OOO,K))=\\
=\varprojlim H^1(F^l(\OOO,K)/F^{l+n}(\OOO,K))=0.
\end{multline*}

Hence $H^1(F^{l}(\OOO,K))=0$ for $l>0$ and
there are surjections
\[
H^0(F^{l}(\OOO,K)) \longrightarrow H^0(\gr^{l}(\OOO,K)) \longrightarrow 0.
\]
This proves our claim. Therefore, near $P$ a general member 
$H\in |\OOO_X|_C$ is given by equations 
$\alpha(y_1,\dots,y_4)=0$ and $\beta(y_1,\dots,y_4)=0$, 
where $\alpha=y_4^2+y_3^3+y_2^3+(\text{terms of degree $\ge 4$})$ \ 
(recall that $\alpha\not\ni y_1^2y_2,\, y_1y_2^2$),
$\beta\equiv \lambda y_3y_1+ y_4\mod F^4(\OOO,K)$, 
and $\lambda\in \OOO_{\CC^4}$ such that
$\lambda(P) \in \CC$ can be chosen arbitrarily.
Hence we can eliminate $y_4$ and get
\[
(H,P) = \{\gamma(y_1,y_2,y_3) =0 \}/\muu_3(1,1,2) \supset C=y_1\text{-axis}/\muu_3,
\]
where $\gamma$ is a $\muu_3$-invariant convergent power series
such that, for $\sigma=(1,1,2)$, $\gamma_{\sigma=3}=y_2^3$ and the term
$\gamma_{\sigma=6}(y_1,0,y_3)$ is squarefree.
Hence we are done by Computation \ref{comp-a3=cube-a6=squarefree}
below.
\end{de}

\begin{de}{\bf Computation.}
\label{comp-a3=cube-a6=squarefree}
Let $(D,P)$ be a normal surface singularity
\[
(D,P)=\{\gamma=0\}/\muu_3\subset \CC^3/\muu_3(1,1,2),
\]
where $\gamma=\gamma(y_1,y_2,y_3)$ is $\muu_3$-invariant and let 
$C:=(y_1\text{-axis})/\muu_3$.
Let $\sigma$ be the weight $(1,1,2)$. Assume that
$\gamma_{\sigma=3}=y_2^3$ and $\gamma_{\sigma=6}(y_1,0,y_3)$ is squarefree.
Then $D$ has only rational singularities and 
$\Delta(D,C)$ is as follows.
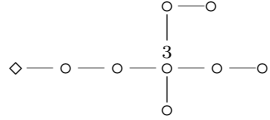
\begin{equation*}\label{comp-a3=cube-a6=squarefree-1}
\begin{array}{lc@{\,}c@{\,}l}
&&\circ&\lin\circ\\
[-4pt]
&&|\\
&\diamond\,\lin\circ\lin\circ
\lin&\stackrel{3}{\circ}&\lin\circ\lin\circ\\
[-4pt]
&&|
\\[-4pt]
&&\circ
\end{array}
\end{equation*}
\end{de}

\begin{proof}
[Sketch of the proof.]
We note that $\gamma_{\sigma=6}(y_1,0,y_3)$ contains $y_3^3$ 
since it is squarefree.
Consider the weighted blowup $\hat H\to H$ with weights $\frac13(1,1,2)$.
The exceptional divisor $\Lambda$ is given by 
$\gamma_{\sigma=3}=y_2^3=0$
in the weighted projective plane 
$\PP(1,1,2)$. Hence $\Lambda$ is a smooth rational curve.
Clearly, $\Sing(\hat H)$ is contained in $\Lambda$.
In the chart $U_1:=\{y_1\neq 0\}$ the surface 
$\hat H$ is given by 
\[
y_2^3+y_1\gamma_{\sigma=6}(1,y_2,y_3)+ y_1^2\gamma_{\sigma=9}(1,y_2,y_3)+\cdots =0.
\]
Hence, $\Sing(\hat H)\cap U_1$ is given by 
$y_1=y_2=\gamma_{\sigma=6}(1,0,y_3)=0$.
Since $\gamma_{\sigma=6}(1,0,y_3)$ is a cubic polynomial 
without multiple factors,
$\Sing(\hat H)\cap U_1$ consists of three points:
$P_0:=(0,0,0)$, $P_1$, $P_2$.
In particular, this shows that $\hat H$ is normal.
Further, 
$\gamma_{\sigma=6}(1,y_2,y_3)$ contains the term $y_3$. Hence
at the origin $\hat H$ has a Du Val singularity of type \typ{A_2}
and the pair 
$$
(\hat H, \Lambda+\hat C)\simeq (\{y_2^3+y_1y_3=0\},\{y_2=0\})
$$
is LC, where $\hat C$ is the proper transform of $C$. This gives us 
the left-hand side of the graph. Similarly, from $P_1$ and $P_2$
we get the upper and the right-hand side of the graph.
The vertex $\circ$ in the bottom comes from the chart $y_3\neq 0$.
The computation of the self-intersection number of the central vertex is 
an easy exercise.
\end{proof}

\begin{de} {\bf Case $\ell(P)=4$ and no type (\typ{III}) points.}
\label{no-type-III-points(4)}
By \cite[2.16]{Mori-1988} $i_P(1)=2$ and (in some coordinate system) $\alpha$
satisfies $\alpha\equiv y_1^4y_3\mod (y_2,y_3,y_4)^2$
(and $C^\sharp$ is the $y_1$-axis).
Hence, $y_2$, $y_4$ form an $\ell$-basis of $\gr_C^1\OOO$.

We prove Claim \ref{(6.3.1)}. Since it has been proved that a type (\typ{III}) point 
does not occur, it remains to settle the $\ell$-isomorphism \ref{(6.3.1.1)}.
If it does not hold, then we have $\gr_C^1\OOO= (2P^\sharp) \toplus (-1)$
and $\gr_C^1\omega = (P^\sharp) \toplus (-2+2P^\sharp)$,
whence $H^1(\gr_C^1\omega) \neq 0$. Thus we get
a contradiction as in \ref{(4)} and Claim \ref{(6.3.1)} is proved.

If $(X,C)$ is flipping, then Claims \ref{(6.3.0)}, \ref{cD3-ip=2-C=ci},
and \ref{(6.3.3)} are already proved in \cite[(6.3)]{Kollar-Mori-1992}.
Since $\ell(P)>2$, $P$ is a double or triple \typ{cD/3} point, and
Claim \ref{converse-cD3-ip=2} is proved in
\cite[(6.3.4)]{Kollar-Mori-1992} if $P$ is a double \typ{cD/3} point.

Assume that $(X,C)$ is not isolated.
Then $P$, as a \typ{cD/3} point, is triple by 
$\ell(P) >2$ and \cite[(6.3.4)]{Kollar-Mori-1992}. 
This proves \ref{(6.3.0)}.

Let $J$ be the $C$-laminal ideal of width 2 such that $J/F_C^2 \OOO
=(0)$ in the $\ell$-splitting \ref{(6.3.1.1)}.
Up to coordinate change we may assume that $y_4$ (resp. $y_2$)
is an $\ell$-free $\ell$-basis of $(0)$ (resp. $(-1+2P^\sharp)$) in
$\gr_C^1 \OOO$ and that $\alpha \equiv y_1^4y_3 \mod
I_C^\sharp J^\sharp$. Whence $y_3 \in F^3(\OOO,J)^\sharp$.
We note that $y_1y_2^2 \not\in\alpha$ in the new coordinates
since $P$ is a triple \typ{cD/3} point.

Since we have $\ell$-isomorphisms
\begin{eqnarray*}
\gr^{2,0}(\OOO,J)&\simeq& \gr^0(\OOO,J)\simeq (0)\\
\gr^{2,1}(\OOO,J) &\simeq& \gr^1(\OOO,J)^{\totimes2} \simeq (-1+P^\sharp),
\end{eqnarray*}
the $\ell$-exact sequence
\[
0 \to \gr^{2,1}(\OOO,J) \to \gr^2(\OOO,J) \to \gr^{2,0}(\OOO,J) \to 0,
\]
is $\ell$-split.
Let $K$ be the ideal such that $J \supset K \supset F^3(\OOO,J)$ and
$K/F^3(\OOO,J)=(0)$ in 
\[
\gr^2(\OOO,J) \simeq (0) \toplus (-1+P^\sharp).
\]
Here we may assume that $y_4$ (resp. $y_2^2$) is an $\ell$-free $\ell$-basis
of $(0)$ (resp $(-1+P^\sharp)$) modulo a coordinate change 
$y_4 \mapsto y_4+(\cdots)y_1y_2^2$. 

We have thus $K^\sharp =(y_2^3,y_3,y_4)$
and 
\[
\gr^1(\OOO,K)=(-1+2P^\sharp), \quad \gr^2(\OOO,K)=(-1+P^\sharp).
\]
We have
$\gr^{3,0}(\OOO,K) \simeq\gr^{2,0}(\OOO,J) \simeq (0)$
and
\[
\alpha \equiv y_1^4y_3 + c y_2^3 \mod I^\sharp K^\sharp,
\]
for some unit $c \in \OOO_X^{\times}$, because
$I^\sharp J^\sharp = I^\sharp K^\sharp+(y_2^3)$ and $y_2^3 \in \alpha$. 
Whence we have an
$\ell$-isomorphism
\[
\gr^{3,1}(\OOO,K) \simeq \gr^1(\OOO,K)^{\totimes 3}(4P^\sharp) \simeq (P^\sharp).
\]
Thus we have an $\ell$-splitting
\[
\gr^3(\OOO,K) \simeq\gr^{3,0}(\OOO,K) \toplus \gr^{3,1}(\OOO,K) \simeq (0) \toplus (P^\sharp).
\]
By a change of coordinate $y_4 \mapsto y_4+ (\cdots)y_1y_3$, 
we may further assume that
$y_4$ (resp. $y_3$) is an $\ell$-free $\ell$-basis of
$(0)$ (resp. $(P^\sharp)$). By the same computation as in 
\ref{(4)}, 
we get
the configuration \ref{(6.2.3.1-b)}.
This contracts to a Du Val point of type \typ{E_6},
and hence $f$ is a divisorial contraction, which proves \ref{(6.3.a)}.
\end{de}

Finally we note that \cite[(6.15) and (6.20)]{Kollar-Mori-1992} 
settled \ref{cD3-ip=2-C=ci} for isolated $(X,C)$,
and \ref{converse-cD3-ip=2} for a double $cD/3$ point. 
We omit the proofs of \ref{cD3-ip=2-C=ci} and \ref{converse-cD3-ip=2}
in other cases since since the arguments are similar.
This completes our treatment of the case $\ell(P)>2$.
\end{proof}

\begin{de}{\bf Examples.}
\label{examples-cD/3}
To show that all the possibilities \ref{eq-theorem-main-cD/3-A2},
\ref{eq-theorem-main-cD/3-D4}, and \ref{eq-theorem-main-cD/3-E6} occur we use deformation arguments.
Consider the surface contraction $f_H: H\to T$ with dual graph 
of the form \ref{eq-theorem-main-cD/3-A2} or
\ref{eq-theorem-main-cD/3-D4}.
By \cite[Proposition (11.4)]{Kollar-Mori-1992} the natural map 
from the deformation space of $H$ to the product of deformation spaces of
singularities $P,\, R\in H$ is smooth, in particular, surjective.
Moreover, the total deformation space $\mathfrak X$ of $H$ has a morphism 
$\mathfrak f$ to the total deformation space $\mathfrak X_Z$ of $T$ so that 
$\mathfrak f|_H=f_H$.
This means in particular 
that any $\QQ$-Gorenstein deformation of singularities of $H$ can be globalized.
Now assume that $(H,P)$ and $(H,R)$
can be obtained as hyperplane sections of 
some terminal singularities $(X,P)$ and $(X,R)$ respectively.
Regard $(X,P)$ and $(X,R)$ as deformation spaces of 
$(H,P)$ and $(H,R)$ respectively. 
By the above there is a globalization $f: X\supset H\to Z\supset T$.

\begin{sde}{\bf Example.}
\label{example-1.1.1}
Consider the surface contraction $f_H: H\to T$ with dual graph 
\ref{eq-theorem-main-cD/3-A2} and consider the following terminal singularities:
\[
\begin{array}{lll}
(X,P)&=&\{y_4^2+y_3^3+y_1y_2(y_1+y_2)=0\}/\muu_3(1,1,2,0),
\\[4pt]
(X,R)&=&\{z_1z_2+z_3^2+z_4^m=0\}, \qquad m\ge 1.
\end{array}
\]
Let $H\subset (X,P)$ is given by $y_4=0$ and $H\subset (X,R)$ 
is given by $z_4=0$. 
By \cite[(6.7.1)]{Kollar-Mori-1992} the dual graph of the 
minimal resolution of $(H,P)$ is the same as that in \ref{eq-theorem-main-cD/3-A2}.
By \ref{examples-cD/3} one obtains the corresponding 
birational contraction $f: X\supset H\to Z\supset T$.
Here $(X,P)$ is a simple \typ{cD/3}-singularity (see 
\cite{Reid-YPG1987}).
Therefore, this $f$ is a divisorial contraction of type \ref{eq-theorem-main-cD/3-A2}.
The point $R\in X$ is smooth if $m=1$ and is a \typ{cA_1}-singularity if $m>1$. 
\end{sde}

\begin{sde}{\bf Example.}
\label{example-1.1.2}
Similarly to Example \ref{example-1.1.1}, take 
\[
(X,P)=\{y_4^2+y_1^2y_2+y_2^6+y_3^3=0\}/\muu_3(1,1,2,0).
\]
By \cite[(6.7.2)]{Kollar-Mori-1992} we get an example of a divisorial contraction
as in \ref{eq-theorem-main-cD/3-D4}.
\end{sde}

\begin{sde}{\bf Example.}
\label{example-1.1.3}
As above, take
\[
(X, P) = \{ y_2^3+y_3^3+y_3y_1^4+y_4^2 \}/\muu_3(1,1,2,0),
\]
where $H$ is cut out by $y_4=y_1y_3$. 
We get an example of a divisorial contraction
as in \ref{eq-theorem-main-cD/3-E6}.
\end{sde}
\end{de}

\section{Case: $P$ is of type \typ{cA/m} and $H$ is normal}
In this section we prove Theorems \ref{theorem-H-not-normal-conic-bundle} and
\ref{theorem-main-birational} in the case where a general 
$H\in |\OOO_X|_C$ is normal. Thus 
throughout this section we assume that $(X,C)$ is an extremal curve germ 
of type \type{IA} or \type{IA^\vee} such that the only non-Gorenstein 
point $P\in X$ is of type \typ{cA/m}
(see \ref{explanations-cA/m-cD/3}).
Let $F\in |{-}K_X|$ be a general member.
Take $H\in |\OOO_X|_C$ so that the pair $(X,F+H)$ is LC
(see Proposition \ref{lemma-lc-near-F}).
Assume that $H$ is normal.
Let $f: (X,C)\to (Z,o)$ be the corresponding contraction.

\begin{thm} {\bf Proposition.}
\label{proposition-H-normal}
In the above notation, 
$H$ has only log terminal singularities of type \typ{T}.
Furthermore, the pair $(H,C)$ is PLT outside of $P$
and $H\setminus \{P\}$ has at most one singular point, 
which if exists is Du Val of type \typ{A_n}. 
If moreover $f$ is birational, then $\Delta(H,C)$ is as in 
\xref{Esq-theorem-H-normal-birational}. 
If moreover $f$ is a $\QQ$-conic bundle, 
then $\Delta(H,C)$ is of the form 
\begin{equation}
\label{eq-semi-stable-theorem}
\circ\lin\circ\lin\circ\lin \bullet\lin \stackrel{4}{\circ}
\end{equation}
In particular, $m=2$ and $(X,P)$ is either a cyclic quotient singularity
$\frac12(1,1,1)$ or a singularity of the form $\{xy+z^2+t^k=0\}/\muu_2(1,1,1,0)$.
\end{thm}

\begin{proof}
First we claim that $H$ has only log terminal singularities.
Write $K_{H}+F|_{H}=(K_X+H+F)|_H\sim 0$. Recall that $F\cap C=\{P\}$. 
So $(H, F|_H)$ is not klt at $P$ and klt at a general point 
of $C$. We see that $(H, F|_H)$ is klt
outside of $P$ by the Connectedness Lemma \cite[5.7]{Shokurov-1992-e-ba},
\cite[17.4]{Utah} (if $f$ is birational) and by Lemma \ref{lemma-connectedness}
(if $f$ is a $\QQ$-conic bundle).
On the other hand, by our assumptions 
and the adjunction formula the pair $(H,F|_H)$ is LC near $F\cap H$,
so the surface $H$ has at worst log terminal singularities.
Further, since $H$ is a Cartier divisor in $X$, the singularities of $H$ are
of type \typ{T} (see \ref{T-singularities-definition}).

Now we claim that the pair $(H,C)$ is PLT outside of $P$.
Assume that $K_H+C$ is not PLT at some point $Q\neq P$.
Take $c$ so that $(H,F|_H+cC)$ is maximally LC.
By the Connectedness Lemma and Lemma \ref{lemma-connectedness} 
we have $c=1$, so $(H,F|_H+C)$ is LC.
Therefore, $H$ has a log terminal singularity at $Q$ and the point 
$(H,Q)$ is Du Val. From the classification of 
log canonical pairs (see, e.g., \cite[Ch. 3]{Utah}) we obtain that 
the part of the dual graph $\Delta(H,C)$ which represents 
$H$ near the singularity $Q$ is of the form
 
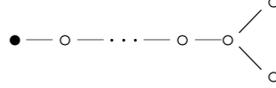
\[
\begin{array}{c@{\,}c@{}c@{\,}c@{\,}}
&&&\circ
\\[-5pt]
&&\diagup
\\[-5pt]
\bullet\,\lin\circ\lin \cdots \lin\circ \lin 
\circ&
\\[-5pt]
&&\diagdown
\\[-5pt]
&&&\circ
\end{array}
\]
But then the corresponding matrix of this subgraph 
is not negative definite, a contradiction.
Thus $(H,C)$ is PLT outside of $P$.
Since any point $Q\in H\setminus \{P\}$ is Gorenstein,
it is Du Val of type \typ{A_n} or smooth.
Near each such a point 
the dual graph $\Delta(H,C)$ is of the form
 \[
\begin{array}{c@{\,}c@{}c@{\,}c@{\,}}
\bullet\,\lin\circ\lin \cdots \lin\circ 
\end{array}
\]
If $(H,C)$ contains two such points, we get a contradiction with 
negative definiteness of the corresponding matrix.
Thus we obtain \ref{Esq-theorem-H-normal-birational}.

Now consider the case where $f$ is a $\QQ$-conic bundle.
If $(H,C)$ is PLT also at $P$, then $H$ has two singularities of types
$\frac1n(1,q)$ and $\frac1n(1, n-q)$ (see Lemma \ref{lemma-plt-log-conic-bundle}).
Since they are of type \typ{T}, we see the following by Proposition 
\ref{T-singularities-characterization}:
\[
(q+1)^2\equiv 0\mod n, \qquad (n-q+1)^2\equiv 0\mod n.
\]
This gives us $4\equiv 0\mod n$. Since 
$X$ is not Gorenstein, the singularities of $H$
are worse than Du Val. Hence, $n=4$. We get the graph \eqref{eq-semi-stable-theorem}.

Finally assume that $(H,C)$ is not PLT at $P$.
Then $\Delta(H,C)$
is of the form \ref{Esq-theorem-H-normal-birational}
with 
$r\neq 1$, $r\neq n$, and $c_1c_n\ge 6$ by Proposition \ref{T-singularities-algorithm}. 
Contracting black vertices
successively, on some step we get a subgraph
\begin{equation}
\label{eq-graph-par-line}
\begin{array}{c@{\,}c@{\,}c@{\,}c@{\,}c@{\,}c@{\,}c@{\,}c@{\,}c@{\,}c@{\,}c@{\,}c@{\,}c}
\stackrel{c_1}{\circ}&\lin&\cdots&\lin&
\stackrel{c_{r-1}}{\circ}&\lin&\bullet&\lin&
\stackrel{c_{r+1}}{\circ}&\lin&\cdots&\lin&
\stackrel{c_n}{\circ}
\end{array}
\end{equation}
Hence strings $[c_{r-1},\dots,c_1]$ and $[c_{r+1}, \dots, c_{n}]$
are conjugate.
This contradicts the following claim because $c_1c_n\ge 6$.
\end{proof}
\begin{sthm}{\bf Claim.}
Let $[a_1,\dots,a_r]$ and $[b_1,\dots, b_s]$ are conjugate strings.
If, for some $c\ge 2$, the string of the form
\begin{newequation}
\label{Esq-string-T-singularities}
[a_r,\dots,a_1,c, b_1,\dots, b_s] 
\end{newequation}
is of type \typ{T}, then it is Du Val.
\end{sthm}
\begin{proof}
Assume that the string 
\eqref{Esq-string-T-singularities} is not Du Val. Take it so that $r+s$ is minimal. 
Since $[a_1,\dots,a_r]$ and $[b_1,\dots, b_s]$ are conjugate, either $a_r=2$ or $b_s=2$. 
Assume that $a_r=2$. 
If $r=1$, then $s=1$ and $b_1=2$, which is a contradiction by 
Proposition \ref{T-singularities-algorithm}, (iii).
Hence, $r>1$, $b_s>2$, and
$[a_{r-1},\dots,a_1,c, b_1,\dots, b_{s-1}, b_s-1]$
is again a non-Du Val \typ{T}-string 
(see Proposition \ref{T-singularities-characterization})
and the 
strings $[a_1,\dots,a_{r-1}]$ and $[b_1,\dots, b_{s-1}, b_s-1]$ are conjugate.
This contradicts our minimality assumption.
\end{proof}

Thus (i) of Theorem \xref{th-index=2} exhausts all $\QQ$-conic bundles with normal $H$.
Explicit examples will be given in Sect. \ref{sect-index-2}.

\begin{de}
In the birational case, similarly to \ref{examples-cD/3} 
any $\QQ$-Gorenstein deformations of singular points of $H$ can be 
globalized by \cite[Proposition (11.4)]{Kollar-Mori-1992}.

\begin{sde}{\bf Example.}
Let $[b_1,\dots,b_r]$ be any \typ{T}-string and let $b_l>2$.
Then the configuration
\[
\begin{array}{l@{\,}c@{\,}l@{\,}}
\stackrel{b_1}{\circ}\lin\cdots\lin&\stackrel{b_l}{\circ}&\lin\cdots\lin\stackrel{b_r}{\circ}
\\[-4pt]
&|&
\\[-4pt]
&\bullet&\lin\underbrace{\circ\lin\cdots\lin\circ}_{k}
\end{array}
\]
where $k\le b_l-3$,
determines a surface germ $(H,C)$ which is contracted to $(T,o)$
with the dual graph 
\[
\begin{array}{l@{\,}l@{\,}l@{\,}}
\stackrel{b_1}{\circ}\lin\cdots\lin&\stackrel{b_l-k-1}{\circ}&\lin\cdots\lin\stackrel{b_r}{\circ}
\end{array}
\]
\end{sde}
For example, for $[b_1,\dots,b_r]=[4]$ and $k=0$, this gives Francia's flip
(see \cite[Theorem 4.7]{Kollar-Mori-1992}).
For $[b_1,\dots,b_r]=[3,2,\dots,2,3]$, $l=r$, and $k=1$ this gives
examples of divisorial extremal neighborhood of index two
\cite[4.7.3.1.1]{Kollar-Mori-1992}.
\end{de}

\section{Case: $P$ is of type \typ{cA/m} and $H$ is not normal}
\begin{de} 
\label{assumptions-non-normal}
In this section we prove Theorems \ref{theorem-H-not-normal-conic-bundle} and
\ref{theorem-main-birational} in the case where a general 
$H\in |\OOO_X|_C$ is not normal. Thus 
throughout this section we assume that $(X,C)$ is an extremal curve germ 
of type \type{IA} or \type{IA^\vee}, the only non-Gorenstein 
point $P\in X$ is of type \typ{cA/m}. 
Let $F\in |{-}K_X|$ be a general member. 
Let $H\in |\OOO_X|_C$ be a non-normal member such that the pair $(X,H+F)$
is LC (see Proposition \ref{lemma-lc-near-F}). 
Let $f: (X,C)\to (Z,o)$ be the corresponding contraction.
\end{de}
\begin{de} \textbf{Setup.}
\label{setup-non-normal}
Let $\nu : H'\to H$ be the normalization and let $\mu :\ti H\to H'$
be the minimal resolution. 
Let $C'=\nu^{-1}(C)$ (with reduced structure) and let $\ti C\subset \ti H$
be the proper transform of $C'$. If $C'$ is reducible,
components of $C'$ (resp. $\ti C$) are denoted by 
$C'_i$ (resp. $\ti C_i$). 
Let $\ov H$ be a minimal model over $T$ (so that $\ov H$ is smooth and 
has no $(-1)$-curves
on fibers over $T$).
Thus we have the following diagram:
\begin{equation*}
\xymatrix{
&\ti H\ar[rr]^{\mu}\ar[dl]_{\ti \upsilon}\ar[ddr]^{\upsilon}&&H'\ar[dr]^{\nu}&&
\\
\ov H\ar[drr]^{\ov \upsilon}&&&& H\ar[dll]_{f_H}
\\
&&T&&&
}
\end{equation*}
Let $\Upsilon:= \nu^{-1}(F\cap H)$.
By \ref{assumptions-non-normal} and Corollary \ref{corollary-structure-H}
we have
\begin{sthm} {\bf Corollary.}
\label{corollary-H-SLC} 
The pair $(H', C'+\Upsilon)$ is LC and 
the restriction map $\nu|_{C'}: C'\to C$ is of degree $2$.
\end{sthm}

\begin{sthm} {\bf Corollary.}
\label{corollary-H-SLC-1} 
Pull-back $C^\sharp$ of $C$ to the index-one cover $(X^\sharp,P^\sharp)\to (X,P)$
is smooth. In particular, $(X,C)$ is of type \type{IA}.
\end{sthm}
\end{de}

Note that $\Delta(H',C')$ is the dual graph of the 1-cycle 
$\upsilon^{-1}(o) \subset \ti{H}$. Hence $\Delta(H',C')$ is
negative semi-definite and its \textit{fundamental cycle} is defined
as usual.

\begin{thm} {\bf Proposition.}
\label{lemma-birational-components-C1}
Under the assumptions of \xref{assumptions-non-normal}
the following are equivalent:
\begin{enumerate}
 \item 
every member of $|\OOO_X|_C$ is 
non-normal,
 \item 
each component of $\ti{C}$ appears with coefficient $>1$ 
in the fundamental cycle $G$ of $\Delta(H',C')$.
\end{enumerate}
In particular, if every member of $|\OOO_X|_C$ is 
non-normal, then all the components of $\ti{C}$ are
contracted by $\ti{\upsilon} : \ti{H} \to \ov H$.
\end{thm}

\begin{proof}
Assume that (ii) does not hold, that is,
a component $\ti{C}_1\subset \ti C$ appears with coefficient 1 in $G$.
Then there is a function $\psi \in \mathfrak{m}_{o, T}$ such that
$\upsilon^*\psi$ has a simple zero along $\ti{C}_1$.
Note that the map $H^0(Z,\OOO_Z) \to H^0(T,\OOO_T)$ is surjective.
Hence $\psi =\phi|_T$ for some $\phi \in \OOO_Z$.
Pick a general point $S \in C$.
If $f^*\phi=0$ is singular along $C$, then $f^*\phi \in I_C^2$
at $S$. 
By commutativity of the above diagram, we have
$\upsilon^*\psi = \mu^*\nu^*(f^*\phi)|_H \in I_{\ti{C}_1}^2$
at a point above $S$. This contradicts the construction
of $\psi$. So $f^*\phi=0$ is smooth along $C$ and a general
member of $|\OOO_X|_C$ is normal, so (i) does not hold. 

Conversely, assume that (i) does not hold.
Then there is a normal member $L\in|\OOO_X|_C$.
Regard $X$ as analytic neighborhood of a general point 
$Q\in C$. Then $H=H_1+H_2$, where $H_1$, $H_2$ are
smooth surfaces intersecting transversely along $C$.
Hence $L$ intersects transversely at least one of $H_1$, $H_2$ along $C$.
This means that $\nu^*L|_H$ is reduced along at least one 
component of $C'$. Thus (ii) does not hold.

As for the last statement, we note that $(T,o)$ is either a 
cyclic quotient singularity (see Proposition \ref{lemma-lc-near-F}) 
or a smooth curve. In both cases $\ti \upsilon (G)$ is reduced.
\end{proof}

\begin{thm} {\bf Proposition.}
\label{proposition-cprime-irreducible-graphs}
Under the assumptions of \xref{assumptions-non-normal},
there are only two possibilities
for the dual graph $\Delta(H',C'+\Upsilon)$:

\begin{sde} \textbf{$C'$ is has two irreducible components: $C'=C_1'+C_2'$.}
\label{case-C-prime-reducible}
\par\medskip
$
\Square\lin
\underbrace{\stackrel{a_r}{\circ}\lin\cdots\lin\stackrel{a_1}{\circ}}_{\varDelta_1}
\lin\diamond\lin
\underbrace{\stackrel{c_1}{\circ}\lin\cdots\lin\stackrel{c_l}{\circ}}_{\varDelta_3}
\lin\diamond\lin
\underbrace{\stackrel{b_1}{\circ}\lin\cdots\lin\stackrel{b_s}{\circ}}_{\varDelta_2}
\lin \Square
$
\end{sde}

\begin{sde} \textbf{$C'$ is irreducible.}
\label{case-C-prime-irreducible}
\par\medskip
$
\Square\lin
\underbrace{\stackrel{a_r}{\circ}\lin\cdots\lin\stackrel{a_1}{\circ}}_{\varDelta_1}
\lin\diamond\lin 
\underbrace{\stackrel{b_1}{\circ}\lin\cdots\lin\stackrel{b_s}{\circ}}_{\varDelta_2}
\lin \Square
$
\end{sde}
\par\noindent
Here $\Square$ corresponds to an irreducible 
component of $\Upsilon$,
$\diamond$ corresponds to an irreducible component of $C'$,
the chain $\varDelta_1$ \textup(resp., $\varDelta_2$\textup) 
corresponds to the singularity of type $\frac1m(1,a)$
\textup(resp., $\frac1m(1,-a)$\textup), 
and in case \xref{case-C-prime-reducible} the chain $\varDelta_3$ 
corresponds to the point $(H',Q')$, where $Q'=C_1'\cap C_2'$.
The strings $[a_1,\dots,a_r]$ and $[b_1,\dots,b_s]$ are conjugate.
If $f$ is birational, then at least one of the vertices 
$\diamond$ corresponds to a $(-1)$-curve under the extra 
assumption that every member of $|\OOO_X|_C$ is non-normal.
If $f$ is a $\QQ$-conic bundle, then all the vertices $\diamond$ correspond 
to $(-1)$-curves.
\end{thm}

\begin{proof}
Note that $C'$ is a fiber of a contraction $H'\to T\ni o$, where 
$(T,o)$ is either a cyclic quotient singularity (see Lemma \ref{lemma-lc-near-F}) 
or a curve germ. 
Hence $p_a(C')=0$ and all components of $C'$ are smooth rational curves.
By Corollary \ref{corollary-H-SLC} $C'$ has at most two components.
So either $C'\simeq \PP^1$ or $C'$ is a union of two $\PP^1$'s 
meeting each other at one point, say $Q'$.

By the classification of log canonical pairs
(see, e.g., \cite[Ch. 3]{Utah})
$\Upsilon$ is smooth at any point $\Upsilon\cap C'$.
On the other hand, $\Upsilon= \nu^{-1}(F\cap H)$, where $H$ is Cartier 
and the pair
$(F,H\cap F)$ is LC.
Hence $\Upsilon$ has exactly two components $\Upsilon_1$, $\Upsilon_2$
and these components are smooth.

Further, since $(H',\Upsilon+C')$ is LC, through any point 
of $H'$ pass at most two 
components of $\Upsilon+C'$. 
Thus for the configuration of $\Upsilon +C'$ on $H'$ 
we have only the following two possibilities:
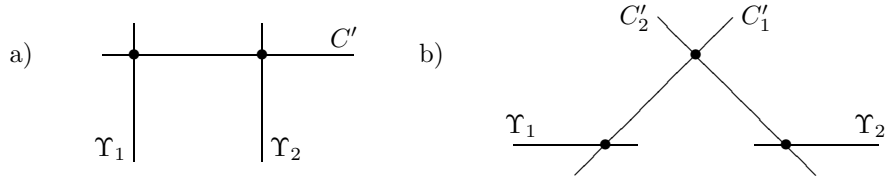
\[\mathrm{a)}\qquad
\xy
(0,0)="Up1n" *{}*+!DR{},
(0,-10)="Up1k" *{}*+!UR{\Upsilon_1},
(17,0)="Up2n" *{}*+!DR{},
(17,-10)="Up2k" *{}*+!UL{\Upsilon_2},
(0,0)="C1n" *{}*+!DL{},
(25,0)="C1k" *{}*+!DL{C'},
{"C1n";"C1k":"Up1n";"Up1k",x} ="I" *{\bullet}*+!UR{},
{"C1n";"C1k":"Up2n";"Up2k",x} ="I" *{\bullet}*+!UR{},
"Up2n";"Up2k"**{} +/1pc/;-/1pc/ **@{-},
"Up1n";"Up1k"**{} +/1pc/;-/1pc/ **@{-},
"C1n";"C1k"**{} +/1pc/;-/1pc/ **@{-}
\endxy
\qquad
\mathrm{b)}\qquad
\xy
(2,2)="C1n" *{}*+!DL{},
(-13,-13)="C1" *{}*+!UL{},
(-2,2)="C2n" *{}*+!DR{},
(13,-13)="C2" *{}*+!UR{},
(-12,-12)="Up1" *{}*+!DR{},
(-20,-12)="Up1k" *{}*+!DR{\Upsilon_1},
(12,-12)="Up2" *{}*+!DR{},
(20,-12)="Up2k" *{}*+!DL{\Upsilon_2},
"C1n";"C1"**{} +/1pc/;-/1pc/ **@{-}*+!L{C_1'},
"C2n";"C2"**{} +/1pc/;-/1pc/ **@{-}*+!R{C_2'},
"Up1";"Up1k"**{} +/1pc/;-/1pc/ **@{-},
"Up2";"Up2k"**{} +/1pc/;-/1pc/ **@{-},
{"Up1";"Up1k":"C1";"C1n",x} ="I" *{\bullet}*+!UR{},
{"Up2";"Up2k":"C2";"C2n",x} ="I" *{\bullet}*+!UR{},
{"C2n";"C2":"C1n";"C1",x} ="I" *{\bullet}*+!U{},
\endxy
\]
Since the pair $(H',\Upsilon+C')$ is LC,
from the classification of log canonical pairs (see, e.g., \cite[Ch. 3]{Utah}) we get 
the desired graphs \ref{case-C-prime-reducible} and \ref{case-C-prime-irreducible}.

It remains to prove the last statements about $(-1)$-curves.
If $f$ is birational, then by Proposition \ref{lemma-birational-components-C1}
at least one of the components of $C'$ is a $(-1)$-curve.
Assume that $f$ is a $\QQ$-conic bundle.
Clearly, the fiber $\upsilon^{-1}(o)$ of a rational curve fibration $\upsilon$ 
contains a $(-1)$-curve and
this curve must coincide with a component of $C'$.
So we are done if $C'$ is irreducible. 
Consider the case \ref{case-C-prime-reducible}.
By the above one of the $\diamond$-vertices corresponds to 
a $(-1)$-curve. Hence the chain 
$\varDelta_1\lin\bullet\lin\varDelta_3\lin\diamond\lin\varDelta_2$
forms a fiber of a rational curve fibration and we may assume that 
$\bullet$ is the only $(-1)$-vertex.
In this case, the chain $\varDelta_1$ is conjugate to both
$\varDelta_2$ and 
$\varDelta_3\lin\diamond\lin\varDelta_2$ 
(see Lemma \ref{lemma-plt-log-conic-bundle}), a contradiction.
\end{proof}

\begin{thm} {\bf Lemma.}
\label{lemma-embedded-dimension}
Let $Q\in H \setminus \{P\}$ 
be any point and let $Q'\in \nu^{-1}(Q)$.
Then 
$4\ge\embdim (H,Q)\ge \embdim (H',Q')-1$.
\end{thm}

\begin{proof}
By Corollary \ref{corollary-H-SLC} the conductor ideal
coincides with the ideal sheaf $I_{C'}$.
The natural map $\OOO_H\to \nu_*\OOO_{H'}$ induces an isomorphism
$I_C \simeq \nu_*I_{C'}$ (any regular function on $H'$ that
vanishes on $C'$ descends to $H$). From the following 
commutative digram
\[
\xymatrix{
0\ar[r]&\nu_*I_{C'}\ar[r]&\nu_*\OOO_{H'}\ar[r]&\nu_*\OOO_{C'}\ar[r]&0
\\
0\ar[r]&I_{C}\ar[r]\ar@{=}[u]&\OOO_{H}\ar[r]
\ar@{^{(}->}[u]
&\OOO_{C}\ar[r]\ar@{^{(}->}[u]&0
}
\]
we have $\nu_*\OOO_{H'}/\OOO_{H}\simeq \nu_*\OOO_{C'}/\OOO_{C}$.
Note that $\nu_*\OOO_{C'}$ is a locally free $\OOO_C$-module
and there is a local splitting $\nu_*\OOO_{C'}=\OOO_C\oplus \OOO_Ct$
for some $t\in \nu_*\OOO_{C'}$. Thus $\nu_*\OOO_{H'}/\OOO_{H}\simeq \OOO_Ct$.
Therefore, $\mathfrak{m}_{Q',H'}/\mathfrak{m}_{Q',H'}^2$ is generated by
$1+\dim \mathfrak{m}_{Q,H}/\mathfrak{m}_{Q,H}^2$ elements
as an $\OOO_{Q,H}$-module.
\end{proof}

\begin{sthm} {\bf Corollary (cf. \cite{Tziolas2005a}).}
\label{corollary-Tziolas2005a}
The chain $\varDelta_3$ in \xref{case-C-prime-reducible}
satisfies the inequality 
\begin{equation}
\label{Esq-lemma-ci}
\embdim (H',Q')-3=\sum (c_i-2)\le 2.
\end{equation}
\end{sthm}
The proof of this statement is contained in 
\cite[Proof of Th. 5.6]{Tziolas2005a}
which is rather computational and uses the classification of degenerate cusp singularities.
Here is a much shorter proof.

\begin{proof}
By Lemma \ref{lemma-embedded-dimension} we have
$\embdim (H',Q')\le \embdim (H,Q)+1\le 5$.
On the other hand, since $(H',Q')$ is a cyclic quotient singularity,
\begin{multline*}
\embdim (H',Q')=-\Bigl(\sum E_i\Bigr)^2+1=
\\
=1+\sum c_i-2 \sum_{i\neq j} E_i\cdot E_j=
3+\sum (c_i-2),
\end{multline*}
where the $E_i$'s are exceptional divisors on the minimal resolution.
This immediately gives the desired inequality.
\end{proof}

\begin{thm} {\bf Proposition.}
\label{lemma-birational-components-C2}
Assume that we are in the case \xref{case-C-prime-reducible}
under \xref{assumptions-non-normal}. Furthermore, assume that
every member of $|\OOO_X|_C$ is non-normal
and $\sum (c_i-2)=2$ \textup(whence $\embdim (H,Q)=4$\textup).
Let $G$ \textup(resp. $G'$\textup) 
be the fundamental cycle of $\Delta(H',C')$ \textup(resp. $\varDelta_3$\textup).
Then $G \ge 2G'$ if and only if $\embdim (M,Q)=4$ for general
$M \in |\OOO_X|_C$.
\end{thm}
\begin{proof}
We have an analytic isomorphism $(H',Q')\simeq \CC^2_{u,v}/\muu_n(1,q)$
for some $n$, $q$ with $\gcd(n,q)=1$.
By Proposition \xref{lemma-birational-components-C1} the graph 
$\ti C_1\lin \varDelta_3\lin \ti C_2$ is contracted on $\ov H$.
Note that $G$ is $\ov\upsilon$-numerically trivial. 
Thus there is a function $\psi\in \OOO_H$ such that
$\mu^*\nu^*\psi=0$ defines $G$ 
near $\mu^{-1}\nu^{-1}(Q)$.
Hence the lifting of $\nu^*\psi$ to $\CC_{u,v}^2$ is given by 
an invariant monomial $\lambda$ multiplied by a unit. 

Since $\sum (c_i-2)=2$, we see
$\embdim (H',Q')=5$ 
and $\embdim (H,Q)=4$ by 
Corollary \ref{corollary-Tziolas2005a} and Lemma 
\xref{lemma-embedded-dimension}, and
$I_{C'} \subset \mathfrak{m'}_{Q',H'}$ is generated by
exactly three invariant monomials in $u,v$ divisible by $uv$.
Thus every minimal generating set of $I_C \subset \mathfrak{m}_{Q,H}$
induces a minimal generating set of $I_{C'} \subset \mathfrak{m'}_{Q',H'}$
(cf. the proof of Lemma \xref{lemma-embedded-dimension}).
This means that $\embdim (M,Q)<4$ for general $M \in |\OOO_X|_C$
iff $\nu^*\psi$ can be a part of a coordinate of $(H',Q')$.
However since the lifting of $\nu^*\psi$ is an invariant monomial (times a unit),
this happens iff $\nu^*\psi$ equals one of the three monomial
generators of $I_{C'}$.

There are only two series of possibilities for $\Delta(H,C)$ near $Q$:

\[
\begin{array}{c@{\,}c@{\,}c@{\,}c@{\,}c}
\diamond&\lin
\underbrace{\circ\lin\cdots\lin\circ}_{a-2}
\lin
&\stackrel{4}{\circ}&\lin 
\underbrace{\circ\lin\cdots\lin\circ}_{b-2}
\lin 
&\diamond
\end{array}\qquad a,\, b\ge 2
\leqno{(*)}
\]
\[
\begin{array}{c@{\,}c@{\,}c@{\,}c@{\,}c@{\,}c@{\,}c}
\diamond&\lin
\underbrace{\circ\cdots\circ}_{a-2}
\lin&\stackrel{3}{\circ}&\lin 
\underbrace{\circ\cdots\circ}_{b-2}
\lin&\stackrel{3}{\circ}&\lin 
\underbrace{\circ\cdots\circ}_{c-2}
\lin &\diamond
\end{array}\qquad a,\, b,\, c\ge 2
\leqno{(**)}
\]
Each monomial in $I_{C'}$ corresponds to
an effective divisor of $\ti H$ with support $\varDelta_3 \cup \tilde{C}$
which is $\mu$-trivial (i.e., numerically trivial along $\varDelta_3$).
The following table gives three such monomials (or divisors) 
$m_A, m_B, m_C$ for each of $(*)$ and $(**)$.
For instance the numbers of the row $m_A$ shows the 
coefficient of the curve corresponding to the vertex 
in the divisor $m_A$.
\begin{center}
\begin{longtable}{c|clcccrc}
(*)&$\diamond$&\multicolumn{2}{c}{$\circ$ \dotfill$\circ$}
&$\stackrel{4}{\circ}$&\multicolumn{2}{c}{$\circ$ \dotfill$\circ$}&$\diamond$
\\[4pt]
$m_A$&$1$&\multicolumn{2}{c}{\dotfill $1$}&$1$&\multicolumn{2}{c}{$3$ \dotfill}&$2b-1$
\\
$m_B$&$2a-1$&\multicolumn{2}{c}{\dotfill $3$}&$1$&\multicolumn{2}{c}{$1$ \dotfill}&$1$
\\
$m_C$&$a$&\multicolumn{2}{c}{\dotfill $2$}&$1$&\multicolumn{2}{c}{$2$ \dotfill}&$b$
\\
\\[-9pt]
\hline
\\[-9pt]
(**)&$\diamond$&
$\circ\ \cdots\ \circ$
&$\stackrel{3}{\circ}$&
$\circ\ \cdots\ \circ$
&$\stackrel{3}{\circ}$&
$\circ\ \cdots\ \circ$
&$\diamond$
\\[4pt]
$m_A$&$1$ &\dotfill & $1$& \dotfill & $b$& \dotfill &$bc +c-1$
\\
$m_B$&$ab+a-1$ & \dotfill & $b$& \dotfill & $1$& \dotfill &$1$
\\
$m_C$&$a$ & \dotfill & $1$& \dotfill & 1& \dotfill &$c$
\end{longtable} \end{center}
It is clear that none of these monomials belong to $\mathfrak{m}_{Q',H'}^2$
because each vanishes to order 1 at one of the vertices with weight
3 or 4. Hence $m_A, m_B, m_C$ are the monomial generators of $I_{Q'}$.
One can also check that the lifting of $\nu^*\psi$ equals one of $m_A, m_B, m_C$ 
iff one of the vertices of weight 3 or 4 appears with coefficient 1 in $G$ iff
$G \not\ge 2G'$.
\end{proof}

\begin{thm} {\bf Proposition.}
\label{proposition-conic-bundle}
Assume that $f$ is a $\QQ$-conic bundle 
germ such that
every member of $|\OOO_X|_C$ is non-normal. Assume furthermore that 
$H\in|\OOO_X|_C$ is taken to be general. 
Then $C'$ is irreducible. 
\end{thm}

\begin{sde}\textbf{ Remark.}
\label{remark-index-2}
If in the above assumptions $X$ is of index $2$, then
$\Delta(H',C'+\Upsilon)$ is of the form 
\[
\Square\lin\circ\lin\bullet\lin\circ\lin\Square.
\]
\end{sde}

\begin{proof}[Proof of \xref{proposition-conic-bundle}]
Assume that $C'$ is reducible.
Then the dual graph $\Delta(H',C')$ is of the form \xref{case-C-prime-reducible}
with $\diamond^2=-1$.
Clearly the chains $\varDelta_1$ and $\varDelta_2$ are not empty
(otherwise $X$ is Gorenstein).
Since the matrix corresponding to 
$\bullet\lin\varDelta_3\lin\bullet$ is negative definite,
the subgraph $\varDelta_3$ is not Du Val. 
We will use the inequality \eqref{Esq-lemma-ci}.
\begin{sde}
Assume that $r=s=1$. 
Then $a_1=b_1=2$ and the graph
\ref{case-C-prime-reducible} or \ref{case-C-prime-irreducible} is of the form 
\[
\circ\lin\bullet\lin\stackrel{4}{\circ}\lin\bullet\lin\circ
\qquad\text{or}\qquad
\circ\lin\bullet\lin
\stackrel{3}{\circ}\lin\underbrace{\circ\cdots\circ}_{l}\lin\stackrel{3}{\circ}
\lin\bullet\lin\circ\quad l\ge 0.
\]
The fundamental cycle $G$ of $\Delta(H',C')$ is
given by
\[
\stackunder{1}\circ\lin\stackunder{2}\bullet\lin\stackunder{1}{\circ}
\lin\stackunder{2}\bullet\lin\stackunder{1}\circ
\quad\text{or}\quad
\stackunder{1}\circ\lin\stackunder{2}\bullet\lin
\stackunder{1}{\circ}\lin\stackunder{1}\circ\cdots\stackunder{1}\circ
\lin\stackunder{1}{\circ}
\lin\stackunder{2}\bullet\lin\stackunder{1}\circ,
\]
respectively. 
Then by Proposition \ref{lemma-birational-components-C2}
our $H$ is not general enough, a contradiction.

>From now on we assume that $rs>1$.
Since $[a_1,\dots,a_r]$ and $[b_1,\dots,b_s]$
are conjugate, we may assume by symmetry that $a_1=2$, $b_1>2$, and $r>1$. 
\end{sde}

\begin{sde}
Consider the case where the chain $\varDelta_3$ contains exactly one curve with
self-intersection $<-2$. 
Then the graph \xref{case-C-prime-reducible}
has the following form:
\[
\stackrel{a_r}{\circ}\lin\cdots\lin\stackrel{a_1=2}{\circ}
\lin\bullet\lin
\underbrace{\stackrel{}{\circ}\lin\cdots\lin\circ
}_{l_1}
\lin\stackrel{c}{\circ}\lin
\underbrace{\circ\lin\cdots\lin\stackrel{}{\circ}}_{l_2}
\lin\bullet\lin
\stackrel{b_1}{\circ}\lin\cdots\lin\stackrel{b_s}{\circ}
\]
where $c=3$ or $4$. 
Since $a_1=2$, it holds $l_1=0$ because the graph $\circ\lin\bullet\lin\circ$
is not negative definite.
Choose the above configuration so that $c$ is minimal.

If $l_2>0$, then contracting both black vertices we get
\[
\stackrel{a_r}{\circ}\lin\cdots\lin\stackrel{a_2}{\circ}
\lin\bullet
\lin\stackrel{c-1}{\circ}\lin
\underbrace{\circ\lin\cdots\lin\stackrel{}{\circ}}_{l_2-1}
\lin\bullet\lin
\stackrel{b_1-1}{\circ}\lin\cdots\lin\stackrel{b_s}{\circ}
\]
The strings $[a_2,\dots,a_r]$ and $[b_1-1,\dots,b_s]$ at the ends are again conjugate.
This contradicts our minimality assumption because $c'=c-1<4$.

Therefore, $l_1=l_2=0$ and \xref{case-C-prime-reducible} is of the form
\[
\stackrel{a_r}{\circ}\lin\cdots\lin\stackrel{a_1}{\circ}
\lin\bullet\lin
\stackrel{c}{\circ}\lin
\bullet\lin
\stackrel{b_1}{\circ}\lin\cdots\lin\stackrel{b_s}{\circ}
\]
Contracting black vertices we get 
\[
\stackrel{a_r}{\circ}\lin\cdots\lin\stackrel{a_2}{\circ}\lin\bullet
\lin
\stackrel{c-2}{\circ}\lin
\stackrel{b_1-1}{\circ}\lin\cdots\lin\stackrel{b_s}{\circ}
\]
Hence $c=4$ and $a_2\ge 3$.
Again the string $[a_2,\dots,a_r]$ is conjugate 
to both $[b_1-1,\dots,b_s]$ 
and $[c-2,b_1-1,\dots,b_s]$, a contradiction.
\end{sde}

\begin{sde}
Now we consider the case where $\varDelta_3$ contains exactly two $(-3)$-curves. 
Then the graph \xref{case-C-prime-reducible}
has the following form:
\[
\stackrel{a_r}{\circ}\lin\cdots\lin\stackrel{a_1=2}{\circ}
\lin\bullet\lin\stackrel{3}{\circ}\lin
\underbrace{\circ\cdots\circ}_{l_1}
\lin\stackrel{3}{\circ}\lin
\underbrace{\circ\cdots\stackrel{}{\circ}}_{l_2}
\lin\bullet\lin
\stackrel{b_1}{\circ}\lin\cdots\lin\stackrel{b_s}{\circ}
\]
(As above $c_1>2$, since $a_1=2$).
If $l_2>0$, then contracting both black vertices we get
\[
\stackrel{a_r}{\circ}\lin\cdots\lin\stackrel{a_2}{\circ}
\lin\bullet\lin
\underbrace{\circ\cdots\circ}_{l_1+1}
\lin\stackrel{3}{\circ}\lin
\underbrace{\circ\cdots\stackrel{}{\circ}}_{l_2-1}
\lin\bullet\lin
\stackrel{b_1-1}{\circ}\lin\cdots\lin\stackrel{b_s}{\circ}
\]
Here again the strings $[a_2,\dots,a_r]$ and $[b_1-1,\dots,b_s]$ are conjugate.
This contradicts the case considered above.
So, $l_2=0$. Then contracting both black vertices we get
\[
\stackrel{a_r}{\circ}\lin\cdots\lin\stackrel{a_2}{\circ}
\lin\bullet\lin
\underbrace{\circ\cdots\circ}_{l_1+2}
\lin
\stackrel{b_1-1}{\circ}\lin\cdots\lin\stackrel{b_s}{\circ}
\]
As above the string $[a_2,\dots,a_r]$ is conjugate 
to both $[b_1-1,\dots,b_s]$ and $[2,\dots,2,b_1-1,\dots,b_s]$, a contradiction.
\end{sde}
\end{proof}

\begin{thm} {\bf Corollary.}
\label{corollary-conic-H}
Let $f$ be a $\QQ$-conic bundle such that a general member $H\in |\OOO_X|_C$ 
is not normal.
Then the germ $(H,C)$ is analytically isomorphic to the germ 
along the line $L:=\{y=z=0\}$
of the hypersurface given by the following weighted polynomial of degree $2m$
in variables $x$, $y$, $z$, $u$\textup:
\[
\phi:= x^{2m-2a}y^2+x^{2a}z^2+yzu.
\]
in 
$\PP(1,a,m-a,m)$, for some integers $a$, $m$ such that $0<a<m$ and $\gcd(a,m)=1$.
\end{thm}

\begin{proof}
By Proposition \ref{proposition-conic-bundle} $(H,C)$ is 
of type \xref{case-C-prime-irreducible}. Then
it is easy to see that the pair $(H,C)$
up to analytic isomorphism is uniquely defined by the types of singularities 
$\frac1m(1,a)$ and $\frac1m(1,-a)$. On the other hand, 
the hypersurface $\phi=0$ satisfies the 
conditions of \xref{case-C-prime-irreducible}.
\end{proof}

Note that we are interested only in the \textit{germ} of the 
hypersurface $\{\phi=0\}$ along $L$.

\begin{sde} \textbf{Remark.} 
Since the germ $(\{\phi=0\}, L)$ is analytically isomorphic
to our $(H,C)$, there is a rational curve fibration on 
$(\{\phi=0\}, L)$ whose central fiber is $L$. 
One can check that this fibration is given by 
the rational function
\[
s=\frac{y^{m-a}z^a}{x^{2a(m-a)}},
\] 
which is regular in a neighborhood of $L$ in $H$.
\end{sde}

\begin{sthm} {\bf Lemma.}
\label{lemma-contraction-exists}
Let $(H, C)$ be as in \xref{corollary-conic-H} and 
let $s: H\to T$ be the corresponding rational curve fibration.
Let 
$t: X\to \CC$ be a one-parameter smoothing of $(H,C)$ 
in a $\QQ$-Gorenstein family.
If $X$ has only terminal singularities, then $(X,C)$ is a 
$\QQ$-conic bundle germ.
\end{sthm}

\begin{proof}
Let $V:=s^{-1}(o)$ (with the scheme structure) and
let $Z$ be the component of 
the Hilbert scheme of $X$ containing 
the point $o=[V]$ representing $V$. Let 
$\mathfrak{X}\subset X\times Z$ be the corresponding universal family.
We have the following commutative diagram
\[
\xymatrix@!0{
& V \ar@{^{(}->}[ld]\ar@{=}[rr]
& & W \ar@{^{(}->}[dl]
\\
X \ar@{..>}[rrdd]\ar[dd]^{t}
& & \mathfrak X \ar[dd]^{\pi}\ar[ll]_(0.4){p}
\\
& 
& & 
\\
\CC& & Z \ar[ll]
}
\]
where $W:=\pi^{-1}(o)$.
Both $V$ and $W$ are locally complete intersections. Moreover, 
$I_V/I_V^2\simeq \OOO_V\oplus \OOO_V$ and
$I_W/I_W^2\simeq \OOO_W\oplus \OOO_W$.
Since $H^1(V,(I_V/I_V^2)^\vee)=0$, $Z$ is smooth at $o$ and 
there is a natural isomorphism 
$\CC^2\simeq T_{o,Z}\simeq H^0(V,(I_V/I_V^2)^\vee)$.
On the other hand, $H^0(V,(I_W/I_W^2)^\vee)\simeq T_{o,Z}$
because $W$ is a fiber of $\pi$. Therefore, there is a 
natural isomorphism
$H^0(W,(I_W/I_W^2)^\vee)\simeq H^0(V,(I_V/I_V^2)^\vee)$ and 
the natural map $(I_W/I_W^2)^\vee\to (I_V/I_V^2)^\vee$ is also an isomorphism.
Thus $p$ is an isomorphism in a neighborhood of $W$.
By shrinking $\mathfrak X$ and $X$ we may assume that there is a contraction
$X\to Z$ such that the whole diagram is commutative.
\end{proof}

The existence of a $\QQ$-Gorenstein smoothing follows from 
\cite{Tziolas2009}. However in our particular case we can 
construct it explicitly:

\begin{sthm} {\bf Lemma.}
Let $(H,C)$, $m$, $a$ be as in Corollary \xref{corollary-conic-H}.
For $s=(s_1,\ldots,s_5)\in \CC_s^5$, 
hypersurfaces $H_{s} \subset \PP(1,a,m-a,m)$ given by equation
\[
\phi_{s}:=\phi+ s_1x^{2m-a}y+s_2x^{m-a}uy+s_3x^{2m}+s_4x^{m}u+s_5u^{2}=0
\]
form a miniversal deformation family of the germ $C \subset H$.
\end{sthm}

\begin{proof}
We compute $T^1_{qG}(H)$ from the $\QQ$-Gorenstein smoothing 
$H \subset P:=\PP(1,a,m-a,m)$ (cf. \cite[\S 3]{Tziolas2009}). 
By definition, $T^1_{qG}(H)$ has an $\ell$-structure
and $T^1_{qG}(H)^{\sharp}=T^1_{qG}(H^{\sharp}).$
Furthermore, we get an exact sequence
\[
\HHom_H(\Omega^1_P,\OOO_H) \longrightarrow \HHom_H(\OOO_P(-H),\OOO_H) 
\longrightarrow T^1_{qG}(H) \to 0
\]
of sheaves with $\ell$-structures. So $T^1_{qG}(H)=\OOO_P(2m)/G$, where
$G$ is generated by $\phi$ and its derivatives. A direct computation shows that
$x^{2m-a}y$, $x^{m-a}yu$, $x^{2m}$, $x^mu$, $u^2$ form a $\CC$-basis of the vector
space $T^1_{qG}(H)$;\ $x^{2m-a}y$, $x^{m-a}yu$ generate the torsion part
of $T^1_{qG}(H)$ and $x^{2m}$, $x^mu$, $u^2$ generate $T^1_{qG}(H)/(\mt{torsion}) \simeq
\OOO_P(2m) \otimes \OOO_C \simeq \OOO_{\PP^1}(2)$.
\end{proof}

\begin{sde}{\bf Example.}
\label{example-conic-bundle}
Let $\alpha$, $\beta\in \CC$ are some general constant
and let $X$ be the threefold given in $\PP(1,a,m-a,m)\times \CC_t$ by
\[
\phi+(\alpha x^m-u)(\beta x^m-u)t=0. 
\]
Then the singularities of $X$ along the curve $C:=\{y=z=t=0\}$
consists of a cyclic quotient singularity of type $\frac1m (1,a,m-a)$
at $\{x=y=z=t=0\}$ and two (Gorenstein) ordinary double points
at $\{\alpha x^m-u=y=z=t=0\}$ and $\{\beta x^m-u=y=z=t=0\}$.
The contraction $X\to Z$ exists by Lemma \xref{lemma-contraction-exists}.
\end{sde}

Thus Theorem 
\ref{theorem-H-not-normal-conic-bundle} is proved.
Now assume that $f$ is birational.

\begin{thm} {\bf Lemma (\cite[Th. 5.6, (1a)]{Tziolas2005a}).}
If $f$ is birational, then $C'$ 
is reducible and the dual graph $\Delta(H',C')$
is of the form 
\xref{case-C-prime-reducible}.
\end{thm}
\begin{proof}
Assume that $\Delta(H',C')$
is of the form 
\xref{case-C-prime-irreducible}.
Then the chain of smooth rational curves corresponding to 
the graph
$\varDelta_1\lin \bullet\lin \varDelta_2$
is contracted by $\upsilon$. On the other hand,
$\varDelta_1$ and $\varDelta_2$ are conjugate.
By Lemma \xref{lemma-plt-log-conic-bundle} this configuration 
corresponds to a rational curve fibration, i.e., $\upsilon$ is not birational, 
a contradiction. 
\end{proof}

\begin{de}
\label{degenerate-cusp}
The singularity $(H,Q)$ is a so-called \textit{degenerate cusp}
\cite{Shepherd-Barron1983}.
One can define the \textit{fundamental cycle} $\Gamma$ 
of $(H,Q)$ and attach an invariant $\zeta=-\Gamma^2$ to $(H,Q)$ such that
\par\medskip 
$\zeta=1$ $\Longleftrightarrow$ 
\parbox[t]{4cm}{$(H',Q')$ is a smooth point}
$\Longleftrightarrow$ $(H,Q)\simeq \{y^2=x^3+x^2z^2\}$,
 \par\medskip 
$\zeta=2$ $\Longleftrightarrow$ 
\parbox[t]{4cm}
{$(H',Q')$ is a Du Val point of type \typ{A_n}, $n\ge 1$}
$\Longleftrightarrow$ 
$(H,Q)\simeq \{y^2=x^2z^2+x^{n+3}\}$,
\par\medskip 
$\zeta=3$ $\Longleftrightarrow$ 
\parbox[t]{4cm}{$\sum (c_i-2)=1$} 
$\Longleftrightarrow$ 
\parbox[t]{4.7cm}{$(H,Q)\simeq \{xyz=y^{a+3}+z^{b+3}\}$, $a,\, b\ge 0$,}
\par\medskip 
$\zeta=4$ $\Longleftrightarrow$ 
\parbox[t]{4cm}{$\sum (c_i-2)=2$}
$\Longleftrightarrow$ $\embdim (H,Q)=4$,
\par\medskip \noindent
(see \cite[\S 1]{Shepherd-Barron1983}). 
Then by \cite[Th. 3.1, Prop. 3.4]{Tziolas2009} 
we have

\begin{sthm} {\bf Theorem.}
In the above notations, a one parameter smoothing 
of $(H,C)$ with only terminal singularities exists if and only if 
\[
\ti C_1^2 +\ti C_2^2+1
+ 4\delta_ {\zeta,1}
+ 4\delta_ {\zeta,2} 
+ 3\delta_ {\zeta,3} 
+ 2\delta_ {\zeta,4} 
\ge 0,
\]
where $\delta_{i,j}$ is Kronecker's delta.
\end{sthm}

\begin{sde} {\bf Remark.}
One can see that the last inequality is equivalent to
\[
\ti C_1^2 +\ti C_2^2+5-\sum(c_i-2)\ge 0,
\]
where we put $\sum(c_i-2)=0$ if $\varDelta_3$ is empty.
\end{sde}

\begin{sde} {\bf Example.}
Assume that the configuration in \ref{case-C-prime-reducible}
is of the form
\[
\circ\lin\stackrel{4}\circ\lin\stackrel{4}\diamond\lin\stackrel{}\bullet
\lin\circ\lin\circ\lin\stackrel{3}\circ
\]
Then $(X,C)$ is a divisorial extremal neighborhood.
By Proposition \xref{lemma-birational-components-C1} every member of $|\OOO_X|_C$
is non-normal. 
By \ref {degenerate-cusp} this $H$ is general in $|\OOO_X|_C$.
\end{sde}
\end{de}

\section{On index two $\QQ$-conic bundles}
\label{sect-index-2}
In this section we give examples of index two $\QQ$-conic bundles.
Let $y_1$, $y_2$, $y_3$, $y_4$; $u$, $v$ be as in 
Theorem \ref{th-index=2}, and let $X\subset \PP(1,1,1,2)\times \CC^2$ given by
\begin{eqnarray*}
0&=&\alpha_1 y_1^2+\alpha_2 u^e y_4 + (\beta_2 u+v)y_3^2\\
0&=&\alpha_3 (y_2^2+\beta_1 y_1 y_3)+ \alpha_4 u y_3^2 + v y_4,
\end{eqnarray*}
where $\alpha_1, \ldots, \alpha_4 \in \CC$ are general, $\beta_1, \beta_2 \in \CC$
are either zero or general, and $e=1, 2, 3$.
Furthermore, $C \subset X$ is given by $y_1=y_2=u=v=0$.

By Bertini's Theorem, we see that the singular locus, $\Sigma$, of $X$
is contained in $\{u=v=0\}$.
Hence, $\Sigma \subset \{u=v=y_1=y_2=0\}$ and
using notation $[y:z]:=(0:0:y:z)\times (0,0)$, we see 
\begin{eqnarray*}
\Sigma &=&\left\{[y_3:y_4] 
\ \left| \ 
\operatorname{rank}
\left(
\begin{matrix}
0 & \alpha_2 e u^{e-1}y_4+\beta_2 y_3^2 & y_3^2\\
\beta_1 y_3 & \alpha_4 y_3^2 & y_4
\end{matrix}
\right.
\right) \le 1
\right\} \cup\{[0:1]\}
\\
&=&
\begin{cases}
\{[0:1]\} & \text{if $\beta_1\neq0$},\\
\left\{\left[1:{\pm\sqrt{\alpha_4/\alpha_2}}\right], [0:1]\right\} 
 &\text{if $\beta_1=0$, $\beta_2=0$, and $e=1$}, \\
\{[0:1]\} &\text{if $\beta_1=0$, $\beta_2=0$, and $e>1$},\\
\{[1:\alpha_4/\beta_2], [0:1]\} &\text{if 
$\beta_1=0$, $\beta_2\neq0$, and $e>1$}.\\
\end{cases}
\end{eqnarray*}

At $[0:1]$, 
the singularity
$(X,[0:1])$ is a hyper-quotient:
\[
\{\alpha_1 y_1^2+\alpha_2 u^e + \beta_2 u y_3^2 
-\alpha_3 y_2^2 y_3^2 - \alpha_3 \beta_1 y_1 y_3^3 - \alpha_4 u y_3^4=0\}/
\muu_2(1,1,1,0).
\]
By \cite[Cor. 2.1]{Mori-1985-cla}, we see that $(X,[0:1])$ is a terminal singularity of type\\
\begin{itemize}
 \item 
\parbox{50pt}{$\frac12(1,1,1)$} if $e=1$,
 \item 
\parbox{50pt}{\typ{cAx/2}} if $e=2$ (cf. \cite[Thm. 12, (3)]{Mori-1985-cla}),
 \item 
\parbox{50pt}{\typ{cD/2}} if $e = 3$ and $\beta_2\neq0$ \cite[Thm. 23]{Mori-1985-cla},
 \item 
\parbox{50pt}{\typ{cE/2}} if $e=3$ and $\beta_2=0$ (cf. \cite[Thm. 25]{Mori-1985-cla}).
\end{itemize}

Every other singular point, if any, is easily seen to be an ordinary double
point, in particular a type \type{III} point.

\begin{enumerate}
 \item 
Case $\beta_1\neq0$: In this case we can assume $\beta_1=-1$ by
change of coordinate $y_1 \mapsto -y_1/\beta_1$, and we are in case (i)
of Theorem \ref{th-index=2}. In this case, $[0:1]$ is the only singular point and it can be
of type $\frac12(1,1,1)$, \typ{cAx/2}, \typ{cD/2} or \typ{cE/2} as above.

 \item 
Case $\beta_1=0$: In this case we are in case (ii) of Theorem \ref{th-index=2}. 
The type of singularity of $(X,C)$ in our example is
\begin{itemize}
 \item 
\parbox{97pt}{$\frac12(1,1,1)$+\type{III}+\type{III}} if $\beta_2=0$ and $e=1$,
 \item 
\parbox{97pt}{\typ{cAx/2}+\type{III}} if $\beta_2\neq0$ and $e=2$,
 \item 
\parbox{97pt}{\typ{cAx/2}} if $\beta_2=0$ and $e=2$,
 \item 
\parbox{97pt}{\typ{cD/2}+\type{III}} if $\beta_2\neq0$ and $e=3$, and
 \item 
\parbox{97pt}{\typ{cE/2}} if $\beta=0$ and $e=3$.
\end{itemize}
\end{enumerate}
In particular, we have shown that all types of terminal index two singularities can 
appear on $\QQ$-conic bundles as in Theorem \ref{th-index=2}.

\par\bigskip\noindent
\textbf{Acknowledgements.}
The authors are grateful to Professors A. Fujiki and O. Fujino
for helpful discussions.
The paper was written when the second author was visiting 
Research Institute of Mathematical Sciences 
and Mathematical Department of Kyoto University
in 2009. He would like to thank both 
institutions for the hospitality and 
excellent working environment.
The first author was partially supported by JSPS Grant-in-Aid for
Scientific Research (B)(2), No. 20340005.
The second author was
partially supported by grants 
RFBR, No. 08-01-00395-a and Kyoto University's Global COE.

\def\cprime{$'$} \def\mathbb#1{\mathbf#1}


\begin{thebibliography}{KSB88}

\bibitem[Art68]{Artin1968}
M. Artin,
\newblock {On the solutions of analytic equations}.
\newblock {\em Invent. math.}, 5: 277--291 (1968).

\bibitem[Bin81]{Bingener1981}
J. Bingener.
\newblock On the existence of analytic contractions.
\newblock {\em Invent. Math.}, 64(1): 25--67, 1981.

\bibitem[KM92]{Kollar-Mori-1992}
J. Koll{\'a}r and S. Mori.
\newblock Classification of three-dimensional flips.
\newblock {\em J. Amer. Math. Soc.}, 5(3):533--703, 1992.

\bibitem[Kol92]{Utah}
J. ~Koll{\'a}r, editor.
\newblock {\em Flips and abundance for algebraic threefolds}.
\newblock Soci\'et\'e Math\'ematique de France, Paris, 1992.
\newblock Papers from the Second Summer Seminar on Algebraic Geometry held at
 the University of Utah, Salt Lake City, Utah, August 1991, Ast\'erisque No.
 \textbf{211} (1992).

\bibitem[Fuj79]{Fujiki1978/79}
A. Fujiki.
\newblock Closedness of the {D}ouady spaces of compact {K}\"ahler spaces.
\newblock {\em Publ. Res. Inst. Math. Sci.}, 14(1): 1--52, 1978/79.

\bibitem[Gro68]{Grothendieck1968}
A. Grothendieck.
\newblock {\em Cohomologie locale des faisceaux coh\'erents et th\'eor\`emes de
  {L}efschetz locaux et globaux {$(SGA$} {$2)$}}.
\newblock North-Holland Publishing Co., Amsterdam, 1968.
\newblock Augment{\'e} d'un expos{\'e} par Mich{\`e}le Raynaud, S{\'e}minaire
  de G{\'e}om{\'e}trie Alg{\'e}brique du Bois-Marie, 1962, Advanced Studies in
  Pure Mathematics, Vol. 2.

\bibitem[Kaw88]{Kawamata-1988-crep}
Y. Kawamata.
\newblock Crepant blowing-up of {$3$}-dimensional canonical singularities and
  its application to degenerations of surfaces.
\newblock {\em Ann. of Math. (2)}, 127(1):93--163, 1988.

\bibitem[KSB88]{Kollar-ShB-1988}
J. ~Koll{\'a}r and N.~I. Shepherd-Barron.
\newblock Threefolds and deformations of surface singularities.
\newblock {\em Invent. Math.}, 91(2):299--338, 1988.

\bibitem[LW86]{Looijenga-Wahl-1986}
E. Looijenga and J. Wahl.
\newblock Quadratic functions and smoothing surface singularities.
\newblock {\em Topology}, 25(3):261--291, 1986.

\bibitem[Mor85]{Mori-1985-cla}
S. Mori.
\newblock On {$3$}-dimensional terminal singularities.
\newblock {\em Nagoya Math. J.}, 98:43--66, 1985.

\bibitem[Mor88]{Mori-1988}
S. Mori.
\newblock Flip theorem and the existence of minimal models for {$3$}-folds.
\newblock {\em J. Amer. Math. Soc.}, 1(1):117--253, 1988.

\bibitem[Mor02]{Mori-2002}
S. Mori.
\newblock On semistable extremal neighborhoods.
\newblock In {\em Higher dimensional birational geometry (Kyoto, 1997)},
{\em Adv. Stud. Pure Math.} ~\textbf{35}, 157--184. Math. Soc. Japan,
 Tokyo, 2002.

\bibitem[MP08]{Mori-Prokhorov-2008}
S. Mori and Y. Prokhorov.
\newblock On {$\mathbf Q$}-conic bundles.
\newblock {\em Publ. Res. Inst. Math. Sci.}, 44(2):315--369, 2008.

\bibitem[MP09]{Mori-Prokhorov-2008III}
S.~Mori and Y. Prokhorov.
\newblock On $\mathbf {Q}$-conic bundles, {III}.
\newblock {\em Publ. Res. Inst. Math. Sci.}, 45:787--810, 2009.

\bibitem[Nak87]{Nakayama1987}
N. Nakayama.
\newblock The lower semicontinuity of the plurigenera of complex varieties.
\newblock In {\em Algebraic geometry, Sendai, 1985},   {\em Adv.
  Stud. Pure Math.} \textbf{10}, 551--590. North-Holland, Amsterdam, 1987.

\bibitem[Pro97]{Prokhorov-1997_e}
Y. Prokhorov.
\newblock On the complementability of the canonical divisor for {M}ori
 fibrations on conics.
\newblock {\em Sbornik. Math.}, 188(11):1665--1685, 1997.

\bibitem[Pro01]{Prokhorov-2001}
Y. Prokhorov.
\newblock {\em Lectures on complements on log surfaces}, {\em MSJ
 Memoirs} ~\textbf{10}.
\newblock Mathematical Society of Japan, Tokyo, 2001.

\bibitem[Rei83]{Reid1983}
M. Reid.
\newblock Minimal models of canonical {$3$}-folds.
\newblock In {\em Algebraic varieties and analytic varieties (Tokyo, 1981)},
{\em Adv. Stud. Pure Math.} ~\textbf{1}, 131--180. North-Holland,
  Amsterdam, 1983.

\bibitem[Rei87]{Reid-YPG1987}
M. Reid.
\newblock Young person's guide to canonical singularities.
\newblock In {\em Algebraic geometry, Bowdoin, 1985 (Brunswick, Maine, 1985)},
{\em Proc. Sympos. Pure Math.} ~\textbf{46}, 345--414. Amer. Math.
 Soc., Providence, RI, 1987.

\bibitem[SB83]{Shepherd-Barron1983}
N.~I. Shepherd-Barron.
\newblock Degenerations with numerically effective canonical divisor.
\newblock In {\em The birational geometry of degenerations ({C}ambridge,
 {M}ass., 1981)}, {\em Progr. Math.} ~\textbf{29}, 33--84. Birkh\"auser
 Boston, Boston, MA, 1983.

\bibitem[Sho93]{Shokurov-1992-e-ba}
V.~V. Shokurov.
\newblock 3-fold log flips.
\newblock {\em Russ. Acad. Sci., Izv., Math.}, 40(1):95--202, 1993.

\bibitem[Ste88]{Stevens-1988}
J.~Stevens.
\newblock On canonical singularities as total spaces of deformations.
\newblock {\em Abh. Math. Sem. Univ. Hamburg}, 58:275--283, 1988.

\bibitem[Tzi05]{Tziolas2005a}
N. Tziolas.
\newblock {Three dimensional divisorial extremal neighborhoods.}
\newblock {\em Math. Ann.}, 333(2):315--354, 2005.

\bibitem[Tzi09]{Tziolas2009}
N. Tziolas.
\newblock {Q-{G}orenstein deformations of nonnormal surfaces}.
\newblock {\em Amer. J. Math.}, 131(1):171--193, 2009.

\end{thebibliography}
\end{document}